\documentstyle[amscd]{amsart}
\newcommand{\functor}[1]{\mbox{{\sc #1}}}
\newcommand{\stack}[1]{\functor{#1}}
\newcommand{\splice}{\quad \text{ and } \quad}
\newcommand{\commasplice}{, \quad}
\newcommand{\pn}{\mbox{${\Bbb P}^N$}}
\newcommand{\pone}{\mbox{${\Bbb P}^1$}}
\newcommand{\len}{\mbox{${\text{length}}$}}
\newcommand{\hilb}{\mbox{${\text{Hilb}}^{dz-g}_{\pn}$}}
\newcommand{\mg}{\mbox{${\cal M}_g$}}
\newcommand{\mgbar}{\mbox{$\overline{\mg}$}}
\newcommand{\spin}{\mbox{${\functor{Spin}}_{r,g}$}}
\newcommand{\spinbar}{\mbox{${\overline{\functor{Spin}}}_{r,g}$}}

\newcommand{\qspin}{\mbox{${{\functor{QSpin}}}_{r,g}$}}
\newcommand{\purespin}{\mbox{${{\functor{Pure}}}_{r,g}$}}

\newcommand{\isom}{\mbox{${\functor{Isom}}$}}
\newcommand{\Isom}{\mbox{${\stack{Isom}}$}}
\newcommand{\spec}[1]{\mbox{${\text{Spec}}(#1)$}}
\newcommand{\proj}[1]{\mbox{${\text{Proj}}(#1)$}}
\newcommand{\quot}{\mbox{${\text{Quot}}$}}
\newcommand{\cx}{\mbox{${\cal X}$}}
\newcommand{\cl}{\mbox{${\cal L}$}}
\newcommand{\ce}{\mbox{${\cal E}$}}
\newcommand{\cf}{\mbox{${\cal F}$}}
\newcommand{\co}{\mbox{${\cal O}$}}
\newcommand{\cg}{\mbox{${\cal G}$}}
\newcommand{\cm}{\mbox{${\cal M}$}}
\newcommand{\cn}{\mbox{${\cal N}$}}
\newcommand{\cs}{\mbox{${\cal S}$}}
\newcommand{\crr}{\mbox{${\cal R}$}}
\newcommand{\tbar}{\mbox{${\overline{t}}$}}
\newcommand{\Hom}{\mbox{$\text{Hom}$}}
\newcommand{\shfhom}{\mbox{${\cal H}\kern-.15em om$}}
\newcommand{\tensor}{\otimes}
\newcommand{\cross}{\times}
\newcommand{\Aut}{\mbox{${\stack{Aut}}$}}
\newcommand{\aut}{\mbox{${\functor{Aut}}$}}
\newcommand{\ann}[1]{\mbox{${\text{Ann}}(#1)$}}
\newcommand{\coker}{\mbox{${\text{coker}}$}}
\newcommand{\im}[1]{\mbox{${\text{im}}(#1)$}}
\newcommand{\ext}{\mbox{${\ce\kern-.15em xt}$}}
\newcommand{\Ext}{\mbox{${\text{Ext}}$}}
\newcommand{\rxy}{\mbox{$R[[x,y]]/(xy-\pi)$}}
\newcommand{\maxid}{\frak m}
\newcommand{\zr}{\mbox{${\Bbb Z}[1/r]$}}
\newcommand{\fp}{{\frak p}}

\newcommand{\th}{^{\text{th}}}
\newcommand{\rth}{\mbox{$r\th$}}
\newcommand{\irightarrow}{@>{\sim}>>}
\newtheorem{theorem}{Theorem}[subsection]
\newtheorem{proposition}[theorem]{Proposition}
\newtheorem{lemma}[theorem]{Lemma}
\newtheorem{corollary}[theorem]{Corollary}

\theoremstyle{remark}

\theoremstyle{definition}
\newtheorem{defn}[theorem]{Definition}
\begin{document}
\title[Generalized Spin Curves]{Torsion-Free Sheaves and
Moduli of Generalized Spin Curves}
\author{Tyler J. Jarvis}
\address{Department of Mathematics\\ Mississippi State University,
MS 39762-5921}
\email{jarvis@@math.MsState.edu}
\thanks{This work was supported in part through a National Defense Science
and Engineering Graduate Fellowship from DARPA and through the
Graduate Assistance in Areas of National Need program of the
U.S. Department of Education.}
\maketitle

\section*{Introduction}

This article treats the compactification of the space of higher spin
curves, i.e.  pairs $(X,{L})$ with ${L}$ an $\rth$ root of the
canonical bundle of $X.$  More precisely, for positive integers $r$
and $g,$ with $g > 2,$ $r$ dividing $2g-2,$ and for a flat family of
smooth curves $f:\cx
\rightarrow T,$ an {\em r-spin structure} on \cx\ is a line bundle
\cl\ such that $\cl^{\tensor r} \cong \omega_{\cx/T}.$
 And an {\em r-spin curve} over $T$ is a flat family of smooth curves
with an $r$-spin structure.  Now, for a fixed base scheme $S$ over
\zr, let \spin\ be the sheafification of the functor which takes an
$S$-scheme $T$ to the set of isomorphism classes of $r$-spin curves
over $T.$  A compactification of the space of spin curves is a space
(scheme or algebraic stack), which is proper over \mgbar\ (the
Deligne-Mumford compactification of the space of curves), and whose
fibre over \mg\ represents, at least coarsely, the functor \spin.

It is possible (see \cite{thesis}) to compactify \spin\ using
geometric invariant theory.  Namely, in the style of L. Caporaso
\cite{cap:thesis}, for a fixed $d>\! >0$ one can choose a subscheme of
the Hilbert scheme \hilb\ with a geometric quotient that coarsely
represents \spin.  And using results of Gieseker (c.f.
\cite[Theorems 1.0.0 and 1.0.1]{gieseker}, ) one can show that the
semi-stable closure of the subscheme in \hilb\ has a categorical
quotient that provides a compactification.  This compactification is
actually a subscheme of Caporaso's compactification of the relative
Picard scheme over \mgbar.

The principle drawback to the GIT compactification is that it is not
obviously the solution to a moduli problem, and therefore it is
difficult to describe the resulting space and to make the construction
work over a general base, rather than only over algebraically closed
fields.  Moreover, the GIT construction requires that one make some
arbitrary choices, and it is not clear that the resulting
compactification is completely independent of these choices.
Therefore, the approach we take here is to pose a moduli problem,
using torsion-free sheaves, and then show that the associated stack is
actually algebraic and that it does indeed compactify \spin.

We discuss three different moduli problems that provide
compactifications and describe some of their characteristics.  The
naive approach would be to use a rank-one torsion-free sheaf $\ce$
with a suitable $\co_{\cx}$-module homomorphism from $\ce^{\tensor r}$
to the canonical bundle.  But this doesn't quite work, as the resulting
space is not separated.  Some additional conditions on the cokernel of
the homomorhism are necesary to make the stack separated, and the
resulting moduli space is called the space of {\em quasi-spin} curves.
The moduli of quasi-spin curves is relatively easy to construct, but
is difficult to describe, due to the presence of nilpotent elements.
Two better moduli problems are {\em spin curves} and {\em pure spin
curves}.  These are quasi-spin curves with some additional local
conditions.  The local conditions require the use of log structures,
and thus the construction of the moduli space of spin curves and pure
spin curves is more difficult than for quasi-spin curves.  But the
resulting spaces have well behaved singularities.  In fact, the space
of pure spin curves is smooth over \mgbar.

\subsection{Previous Results}

Work on this problem in the special case where the base $S$ is ${\Bbb
C}$ and $r=2$ has been done by M. Cornalba in \cite{corn:theta} and
over a more general base by P. Deligne in \cite{deligne:letter}.  P.
Sipe and C.J. Earle have studied $\rth$ roots of the canonical
bundle on the universal Teichm\"uller curve (c.f.
\cite{sipe:roots2,sipe:roots} and \cite{sipe:roots3}).  And
topological properties of the uncompactified moduli space of $2$-spin
curves have been studied in many places (e.g. \cite{har:spin2,lebrun}
).

\subsection{Overview}

In the first section we present some background on torsion-free
sheaves, some geometric motivation for their use, and some results of
Faltings on the local structure of torsion-free sheaves.  In the
second section we define the first moduli problem---singular {\em
quasi-spin curves}---and make some local calculations that lead
naturally to another condition we impose later on quasi-spin curves to
get singular {\em spin curves}.  In the third section we discuss how
to move from local to global structures using log-structures, and we
use them to formally define the spin curves and pure-spin curves.  The
fourth and fifth sections treat deformation theory and isomorphisms,
respectively.  The sixth section covers the construction of the
different compactifications, proves that they are algebraic stacks,
and discusses the nature of their singularities.  And in the last
section we prove that all the constructions are proper over \mgbar,
and therefore are true compactifications.

\subsection{Notation and Conventions}

We will use the term {\em semi-stable curve} of genus $g$ to mean a
flat, proper morphism $ X \rightarrow T $ whose geometric fibres $X_t$
are reduced, connected, one-dimensional schemes, with only ordinary
double points, and with $\dim H^1(X_{t}, \co_{X_{t}}) = g.$  A {\em
stable curve} is a semi-stable curve of genus greater than one, with
the additional property that any irreducible component which is
isomorphic to \pone\ meets the rest of the curve in at least three
points.  Irreducible components of a semi-stable curve which are
isomorphic to \pone\ but meet the curve in only two points will be
called {\em exceptional} curves.  By {\em line bundle} we mean an
invertible (locally free of rank one) coherent sheaf.  An $r$-spin
structure on a smooth curve $X/T$ will be a line bundle $\cl$ such
that $\cl^{\tensor r}$ is isomorphic to the canonical bundle
$\omega_{X_{T}}.$  A smooth $r$-spin curve will be a smooth curve
$X/S$ with a spin structure.

\section{Torsion-Free Sheaves}

Compactifications of curve-line bundle pairs using geometric invariant
theory give boundary points that correspond to pairs $(X,L)$ with $X$ a
semi-stable curve having at most one exceptional curve (copy of \pone\
that intersects the remaining curve in at most two points) in each
chain of exceptional curves, and $L$ a line bundle of degree one on
each exceptional curve in $X.$ Contracting all the exceptional curves
makes the underlying curve stable, and the direct image of $L$ is a
torsion-free sheaf; namely, it has no associated primes of height one.
Furthermore the torsion-free sheaves on stable curves don't have the
problem of having infinite automorphism groups that the line bundles
on semi-stable curves have.  It is therefore natural to expect that
torsion-free sheaves will be well-suited to the compactification of
the moduli of spin curves, and this is, in fact, the case.

To begin, we define torsion-free sheaves.
\begin{defn} By {\em relatively torsion-free sheaf} (or just
torsion-free sheaf) on a stable or semi-stable curve $f:\cx
\rightarrow T\ ,$ we mean a coherent sheaf \ce\ of
$\co_{\cx}$-modules, which is of finite presentation and flat over
$T,$ with the additional property that on each fibre $\cx_{t} = \cx
\cross_T \spec{k(t)}$ the induced $\ce_{t}$ has no associated primes
of height one.  Of course, on the open set where $f$ is smooth, a
torsion-free sheaf is locally free.
\end{defn}

Our ultimate goal is to define and describe a notion of $r$-spin
structure for stable curves that corresponds to the previously defined
notion for smooth curves.  Spin structures on a family of stable
curves $\cx @>>> T$ will be pairs $(\ce,b)$ of a relatively
torsion-free sheaf $\ce$ and a morphism $b:\ce^{\tensor r} @>>>
\omega_{\cx/T}$ of $\co_{\cx}$-modules, having certain properties that
we will describe later.  But before we can define spin structures, we
need some general results about relatively torsion-free sheaves.

\subsection{General Properties of Torsion-Free Sheaves}

The following propositions give several basic but important properties of
torsion-free sheaves.
\begin{proposition} \label{exacttf} Given any
family of semi-stable curves $\cx/T$ and an exact sequence $$ 0
\rightarrow \ce' \rightarrow \ce \rightarrow \ce'' \rightarrow 0$$ of
coherent sheaves on $\cx,$
\begin{enumerate}
\item If $\ce''$
is flat over $T$ and of finite presentation, and if $\ce$ is
relatively torsion-free, then $\ce'$ is also relatively torsion-free.
\item If $\ce'$ and $\ce''$ are relatively torsion-free, then $\ce$ is.
\end{enumerate}
\end{proposition}
\begin{pf}
That the sheaves in question are of finite presentation is
straightforward to check.  It is enough to check that the sheaves are
torsion-free on each fibre,  and the flatness of $\ce''$ over $T$ means that
the sequence is still exact after restriction to the fibres, where the
proposition is clear.
\end{pf}

\begin{proposition}
For any invertible sheaf $\cl$ and any relatively torsion-free sheaf
$\ce$ on $\cx/T,$ the sheaves $\ext_{\cx}^i (\ce, \cl)$ are zero for all
$i >0.$
\end{proposition}
\begin{pf}
By \cite[7.3.1.1]{ega3} it is enough to check this on the individual
fibres, i.e. we may assume that $T$ is a field, and it is enough to
check at the stalk of a closed point $\fp.$  But in this case
$$\ext_{X}^{i} (\ce, \cl)_{\fp}  = \Ext^{i}_{\co_{X,\fp}}(\ce_{\fp},
\cl_{\fp})  = \Ext^{i}_{\co_{X,\fp}} (\ce_{\fp}, \co_{\fp}).$$

And $X$ is Gorenstein, so these vanish for all $i >1.$  And in the case
$i =1,$ by duality theory (\cite[Theorem 6.3]{localcoho}), $$
\Ext^{1}_{\co_{\fp}} (\ce_{\fp}, \co_{\fp}) \irightarrow
\Hom_{\co_{x,\fp}} (H^{0}_{\{\fp\}} (\ce_{\fp}), I)$$ for some
dualizing module $I.$ But $\ce_{\fp}$ is torsion-free, so it has no
elements with support equal to $\{\fp\},$ so
$H^{0}_{\{\fp\}}(\ce_{\fp})=0$ and thus also
$\Ext^{1}_{\co_{X,\fp}}(\ce_{\fp}, \co_{x, \fp}) = 0.\qed$
\renewcommand{\qed}{}
\end{pf}

\begin{proposition}
If $\ce$ is relatively torsion-free, then for any invertible sheaf
$\cl$ the sheaf $\shfhom_{\cx} (\ce, \cl)$ is
also relatively torsion free.
\end{proposition}
\begin{pf}
$\shfhom (-, \cl)$ preserves flatness over $R$ and commutes with base
change, so it suffices to check the proposition over a field, where it
is clear.
\end{pf}

\begin{proposition}
  Any relatively torsion-free sheaf $\ce$ is reflexive.
\end{proposition}
\begin{pf}
This follows from local duality.
\end{pf}

\subsection{Torsion-Free Sheaves on Semi-Stable Curves over a Field}

 It is well known that the stalk $\cf_{\fp}$ of a rank-one
torsion-free sheaf $\cf$ at a singular point $\fp$ of a semi-stable
curve $X$ is isomorphic either to $\co_{X,\fp}$ or to the maximal
ideal $\maxid_{\fp},$ which is isomorphic to the direct image $\pi_{*}
\co_{\tilde{X},\fp}$ of the normalization $\co_{\tilde{X},\fp}.$ In
particular, if the completion $\hat{\cf}_{\fp}$ of the stalk at $\fp$
is not free, then $\hat{\cf}_{\fp} \cong xk[[x]]\oplus yk[[y]]$ over
$\hat{\co}_{X,\fp} \cong k[[x,y]]/xy,$ where $k$ is the residue field
$\co_{X,\fp}/\maxid_{\fp}.$

The following are some simple but useful results which describe
torsion-free sheaves in terms of line bundles on the normalization of
the curve $X.$ Namely, if $\pi : X^{\nu} @>>> X $ is the normalization
of $X$ at one point $\bar \fp$ of the singular set of $\cf$ (i.e.
where $\cf$ is not free), then the $\co_{X^{\nu}}$-module $\pi^{*}
\cf$ has torsion elements, but modulo the torsion elements it is free
near the two points of $\pi^{-1}(\overline \fp).$ For any
quasi-coherent sheaf $\cg,$ define $\pi^{\natural} \cg$ to be the
torsion-free $\co_{X^{\nu}}$-module $(\pi^{*}\cg/\text{torsion}).$
Straighforward checking yields the following proposition.

\begin{proposition}\label{equiv}
In the situation above, where $\pi$ is the normalization of $X$ at a
singularity of a rank-one, torsion-free sheaf $\cf,$ the canonical map
$\cf \rightarrow \pi_{*} \pi^{\natural} \cf$ is an isomorphism.
\end{proposition}

As usual, define the degree of a sheaf \cf\ on a curve $Y$ over an
algebraically closed field to be $\deg (\cf) = \chi (\cf) -
\chi(\co_{Y}).$ We are primarily interested in the case of relatively
torsion-free sheaves on stable curves, and in this case the degree of
$\ce$ is locally constant on the fibres since $\ce$ is flat over the
base. It is easy to see that this definition of degree corresponds to
the usual definition of degree if $\cf$ is a line bundle.

\begin{proposition}  If $\cf$ is a rank-one torsion-free sheaf and
 $\pi$ is, as above, the normalization of $X$ at a singular point
 of $\cf,$ then $\deg (\pi^{\natural} \cf) = \deg (\cf) -1.$
\end{proposition}
\begin{pf}
$R^i \pi_* (\cf) =0$ for all $i >0,$ so the Leray spectral sequence
degenerates, and $H^{q}(X, \cf) = H^{q} (Y, f_{*} \cf)$ for all $q.$
Thus $\chi (\cf) = \chi (\pi^{\natural}\cf),$ and $\chi( \pi _{*}
\co _{X^{\nu}}) = \chi (\co_{X^{\nu}}).$ Taking Euler-Poincar\'e
characteristics of the exact sequence $ 0\rightarrow \co_{X}
\rightarrow
\pi_{*} \co_{X^{\nu}} \rightarrow k \rightarrow 0 $
gives $\chi (\pi_{*}\co_{X^{\nu}})= \chi (\co_{X}) + 1, $
and thus $\deg (\cf) = \deg (\pi^{\natural} \cf ) + 1.\qed$
\renewcommand{\qed}{}
\end{pf}

\begin{proposition}
 If $ X^{\nu} @>{p}>> X$ is the normalization of $X$ at all the
singularities of $\cf,$ then $p^{\natural} \cf$ is invertible and
$(p^{\natural} \cf)^{\tensor r} \cong p^{\natural} (\cf^{\tensor r}).$
In fact, $ p^{\natural} \ce_1 \tensor p^{\natural} \ce_{2} \tensor
\dots \tensor p^{\natural} \ce_{n} = p^{\natural} (\ce_{1} \tensor
\dots \tensor \ce_{n})$ for any torsion-free sheaves $\ce_1, \ce_2,
\dots, \ce_n$ with singularities equal to the singularities of $\cf.$
\end{proposition}
 This is, again, straighforward to check.

\begin{proposition}
 $\pi^{\natural}$ is a covariant functor from coherent sheaves on $X$
to torsion-free sheaves on $X^{\nu}.$ And if ${\cal{TORF}}_{\nu}$ is
the category of rank-one torsion-free $\co_X$-modules with
singularities exactly those which are normalized by $\pi,$ and
${\cal{PIC}}_{\nu}$ is the category of invertible
$\co_{X^{\nu}}$-modules, then the categories ${\cal{TORF}}_{\nu}$ and
${\cal{PIC}}_{\nu}$ are equivalent via $\pi^{\natural}$ and $\pi_{*}.$
\end{proposition}
\begin{pf}
 The first part is clear except perhaps the fact that $\pi^{\natural}$
commutes with composition of morphisms.  But the effect of
$\pi^{\natural}$ applied to a morphism $ (\ce @>f>> \cf)$ is induced
by the composition $ \pi^{*} \ce @>\pi^{*}f>> \pi^{*} \cf \rightarrow
\pi^{\natural} \cf,$ which factors through $\pi^{\natural} \ce,$ since
the target has no torsion.  Moreover, $ \pi^{\natural}(f \circ g) $ is
given by the following commutative diagram.
\begin{center}
\setlength{\unitlength}{0.012500in}%
\begin{picture}(182,100)(80,677)
\put( 80,760){\makebox(0,0)[lb]{\smash{$\pi^{*}\ce$}}}
\put(160,760){\makebox(0,0)[lb]{\smash{$\pi^{*}\cf$}}}
\put(240,760){\makebox(0,0)[lb]{\smash{$\pi^{*}\cg$}}}
\thicklines
\put(100,765){\vector( 1, 0){ 55}}
\put(180,765){\vector( 1, 0){ 55}}
\put( 90,760){\vector( 0,-1){ 45}}
\put(170,760){\vector( 0,-1){ 45}}
\put(250,760){\vector( 0,-1){ 45}}
\put(100,705){\vector( 1, 0){ 55}}
\put(180,705){\vector( 1, 0){ 55}}
\put(100,699){\vector( 1,0){ 135}}
\put(160,690){\makebox(0,0)[lb]{\smash{$\scriptscriptstyle\pi^{\natural}
 (f\circ g)$}}}
\put(115,710){\makebox(0,0)[lb]{\smash{$\scriptscriptstyle\pi^{\natural} f$}}}
\put(115,770){\makebox(0,0)[lb]{\smash{$\scriptscriptstyle\pi^{*} f$}}}
\put(195,770){\makebox(0,0)[lb]{\smash{$\scriptscriptstyle\pi^{*} g$}}}
\put(195,710){\makebox(0,0)[lb]{\smash{$\scriptscriptstyle\pi^{\natural} g$}}}
\put( 80,702){\makebox(0,0)[lb]{\smash{$\pi^{\natural}\ce$}}}
\put(160,702){\makebox(0,0)[lb]{\smash{$\pi^{\natural}\cf$}}}
\put(240,702){\makebox(0,0)[lb]{\smash{$\pi^{\natural}\cg$}}}
\end{picture}
\end{center}
 Both sets of bottom arrows $ \pi^{\natural} (f \circ g)$ and
$\pi^{\natural} f \circ \pi^{\natural} g$ commute in the rectangle,
and the map $ \pi^{*} \rightarrow \pi^{\natural} $ is surjective, so
the bottom arrows must commute.  The second part of the proposition is
clear by Proposition~\ref{equiv}.
\end{pf}

\subsection{Boundedness of Rank-One Torsion-Free Sheaves}

Since $X$ is projective, we can choose a very ample $ \co(1) $ on $X.$
In general, we will write $\cf(m)$ for $\cf \tensor \co(1)^{\tensor
m}.$ One of the more important facts about torsion-free sheaves of
rank one is that they form bounded families.  In other words, the
following proposition holds.  This will later be important in proving
that the functor of spin curves has a versal deformation.
\begin{proposition}
\label{thm:dsousa}
If \cf is a rank-one torsion-free sheaf on $X,$ there is an integer
$m_0$ depending only on the degree of \cf\ on each irreducible
component of $X,$ and on the genus of $X,$
such that for $m \geq m_0$ the following holds.  \begin{enumerate}
\item $H^1(X,\cf (m)) = 0.$ \item $\cf (m)$ is generated by global
sections. \end{enumerate}
\end{proposition}
This is a straightforward generalization of D'Souza's propositions in
Section Three of \cite{dsou:jac}.

\subsection{Some Geometry}\label{geom}

The original motivation for studying torsion-free sheaves comes from
the fact that boundary points in the GIT compactification correspond
approximately to pairs $(X,\cl)$ with $X$ a semi-stable curve, having
no more than one exceptional curve in each chain of exceptional curves
in any fibre, and $\cl$ is a line bundle with degree one on each of
the exceptional curves (c.f. \cite{cap:thesis,gieseker,thesis}).
Given a curve-bundle pair $(X,\cl)$ of this sort, contracting all of
the exceptional curves, i.e.  projection from $X$ to its stable model
$\rho: X @>>>\bar X$ makes $\rho_*\cl$ into a torsion-free sheaf.  This
contraction and a related ``inverse'' action of blowing up certain
ideals warrant a more careful look.

Let $X/B$ be a stable curve over $B=\spec R,$ with $R$ a complete
local ring.  The completion $A=\hat{\co}_{X,\fp} $ of the local ring
of $X$ at a point $\fp$ in the special fibre is of the form $A=\rxy$
for some choice of $\pi$ in the ring $R.$ If $\fp$ is singular in the
special fibre, then $\pi$ is in the maximal ideal $\maxid$ of $R.$
Given two elements $p$ and $q$ of $R$ such that $pq=\pi,$ we will
construct a semi-stable curve with the properties mentioned above,
namely with no chains of exceptional curves of length greater than
one.  Essentially, we want to undo the contraction to associate a
suitable semi-stable curve and line-bundle to each stable curve and
rank-one torsion-free sheaf.

First notice that if $p$ is not a zero divisor, blowing up the ideal
$I= (x,p)$ in $A$ to get $\tilde{X}_I :=
\proj_{X}(\oplus_{n}I^n)$ gives $$\tilde{X}_I \cong
\proj_{A}(A[P,\Xi]/(Px-p\Xi,Pq-\Xi y)).$$ Similarly, if $q$ is
not a zero divisor, blowing up the ideal $J= (y,q)$ gives
$$\tilde{X}_J \cong \proj_{A}(A [Q,Y]/(Qy-yQ, pQ-xY)).$$
And these are isomorphic via $P \mapsto Y,\ \Xi \mapsto Q.$
In general, define $\tilde{X} (p,q)$ to be the $X$-scheme
defined locally as $$\tilde{X} := \tilde{X} (p,q):=
\proj_{A}(A[P,\Xi]/(Px-p\Xi, Pq-\Xi y)) @>\rho>> X,$$
regardless of whether or not $p$ or $q$ is a zero divisor.  It is
straighforward to check that the curve $\tilde{X}(p,q)$ is actually
semi-stable of the desired form (only one exceptional curve) and has
stable model equal to $X.$

Let $s=\frac{\Xi}{P},$ and let $U$ in $\tilde{X}$ be the open set
$U=\spec{A[s]/(x-ps, ys-q)}.$ Similarly, setting
$t=1/s=\frac{P}{\Xi}$ let $V$ in $\tilde{X}$ be the open set $V=
\spec{A[t]/(xt-p, y-qt)}.$ And finally, let $\tilde{A}$ be the
ring $ \tilde{A} = A[P,\Xi]/(p\Xi-xP,qP-\Xi y).$ The union of $U$ and
$V$ is all of $\tilde{X},$ and there are canonical line bundles on
$\tilde{X},$ namely $\co_{\tilde{X}}(n) $ for all $n$ where
$$\co_{\tilde{X}} (n) = \begin{cases} {P}^{n} \co_U
\text{ on } U \\ {\Xi}^n \co_V \text{ on } V
\end{cases}.$$

\begin{proposition}\label{tilden}
In the above construction of $\tilde{X}$ and $\cl:=
\co_{\tilde{X}}(n)$ the following hold.
\begin{enumerate}
\item $n \geq -1 $ implies that $ \rho_{*}  \cl $ is flat over
      $R,$ it commutes with base change, and $ R^{1} \rho_{*} \cl =0.$
\item For $n \geq 0,$ $\Gamma (\tilde{X}, \cl) \cong \tilde{A}_n,$
      the $n^{th}$ graded piece of $\tilde{A},$ and the natural map
      $ \rho^{*} \rho_{*} \cl \rightarrow \cl$ is surjective.
\item $\rho_{*} \co_{\tilde{X}} = \co_X.$
\item $n=1$ implies $ \rho_{*} \cl$ is torsion-free of rank one.
\end{enumerate}
\end{proposition}
\begin{pf}
 It suffices to consider the case $A = R[x,y]/(xy-\pi),$ and in this
case $\co_U \cong R[s,y]/(sy-q),$ and $\co_V \cong R[t,x]/(xt-p),$
both of which are flat over $R,$ and $\co_{U \cap V}
\cong R[s,t]/(st-1).$

Now $\Gamma (\tilde{X}, \co_{\tilde{X}}(n)) = \{(g,f) \in \co_U \oplus
\co_V | g = s^n f$ on $ \co_{U \cap V} \}$ So $f$ and $g$ are of
the form $$ g = sg_+(s) + g_0 + yg_- (y) \splice f = tf_- (t) + f_0 +
xf_+ (x)$$ with $g_+(s) \in R[s], g_- (y) \in R[y], f_- (t) \in R[t],$
and $f_+(x) \in R[x],$ and $$ sg_+ (s) + g_0 + tqg_- (tq) = s^n (tf_-
(t) + f_0 + spf_+ (sp)).$$ So $$sg_+ (s) + g_0 + tqg_- (tq) = s^{n-1}
f_- (t) + s^n f_0 + s^{n+1} pf_+ (sp),$$ and rewriting $sg_+ (s)$ as $
s^{n+1} \gamma (s) + s^n g_n + \dots + sg_1,$ with $g_i \in R,$ and
$\gamma (s) \in R[s],$ and $tf_+(t)$ as $t^{n+1} \phi (t) + t^{n}
f_{n} + t^{n-1} f_{n-1} + \dots + tf_1,$ with $\phi (t) \in R[t],$
gives $$ s^{n+1} \gamma (s)+ s^n g_n + \dots + sg_1 + g_0 + tqg_- (tq)
= t\phi (t) + f_n + \dots + s^{n-1} f_1 + s^n f_0 + s^{n+1} pf_+ (sp).$$

 Now $s^{n+1} \gamma (s) = s^{n+1} pf_+ (sp),$ and even if $p$ or $q$
is a zero divisor, this implies that $\gamma(s) = pf_+ (sp),$ because
$\gamma(s)$ is an element of $R[s] \subseteq \co_U $ and $\ann
(s^{n+1}) \cap R[s] = (0)$ (i.e. $R[s] \subseteq \co_U \rightarrow
\co_{U \cap V}$ is injective).  Similarly, $g_{n-1} = f_i $ for $
0 \leq i \leq n$ and $\phi (t) = qg_- (tq).$

Several things are easy to see from this
formulation, namely
\begin{enumerate}
\item $\Gamma (\tilde X, \cl)$ is free over $R,$ and hence
$\rho_{*} \cl$ is $R$-flat as long as $n \geq -1.$
\item $\rho_{*}\co_{\tilde{X}} = \co_{X}$
\item If $(g,f)$ is an element of
$\Gamma(\tilde{X}, \cl),$ then $g$ can be written as
$$g = s^n xf_+ (x) + s^n f_0 + s^{n-1} f_1 + \dots + sf_{n-1} + f_n +
yg_- (y)$$ in $\co_U$ (with $ x=sp$), and $f$ can be written as $$f =
xf_+ (x) + f_0 + tf_1 + t^2 f_2 + \dots + t^n f_n + t^n yg_- (y)$$ in
$\co_V$ (with $y = tq$).  \end{enumerate}

 The third fact shows that the element $$(f_0 + xf_+ (x))
{\Xi}^n + f_1 {\Xi}^{n-1} P + \dots + f_{n-1} \Xi {P}^{n-1} + (f_n +
yg_- (y)) {P}^n$$ in $\tilde{A}_n$ maps to $(g,f)$ in
$\Gamma(\tilde{X},\cl),$ and the natural homomorphism $\tilde{A}_n
\rightarrow \Gamma (\tilde{X}, \cl)$ is surjective.
Moreover, it is easy to check that if we write $\tilde{A}_n$ as
$\tilde{A}_n = \{ (F_0 {\Xi}^n + F_1 {\Xi}^{n-1} P \dots + F_n
{P}^n)|F_i \in A\},$ we can assume that $F_0$ is in $R[x]$ and $F_n$
is in $R[y]$ and all the remaining terms are in $R$ (i.e. $F_i \in R$
for $0 < i < n$).  And thus the homomorphism $\tilde{A}_n \rightarrow
\Gamma(\tilde{X}, \cl)$ must actually be injective, hence an isomorphism.

Now since $R^1\rho_{*}$ is right exact, to see that it vanishes, it
suffices to check that it vanishes for each fibre. And $H^1 (\tilde{X}
\cross_{X} x, \co(n))$ is zero for all $x$ in $X,$ except possibly the
singular points of $X.$ But over a singular point  $$H^1
(\co(n)) \cong H^1 (\pone, \co_{\pone} (n)) = 0.$$ And thus $R^1
\rho_{*}$ is zero, and $\rho_{*} \cl$ commutes with base change if $n
\geq -1.$

To show that $ \rho^{*} \rho_{*} \cl \rightarrow \cl$ is surjective,
note that this map is locally (on $U$) just the map taking $ \{ (g,f)
| g = s^{n} f\}
\tensor_{A} \co_U$ to $ \co_U,$ given by $(g,f) \tensor z \mapsto gz.$
So it suffices to show that there exists $(g,f) \in \rho_{*} \cl$ such
that $ g =1.$ But $ (g,f) = (1,t^{n})$ works as long as $t^{n} \in
\co_V,$ i.e. if $ n \geq 0.$ A similar computation holds over $V,$ so
$\rho^{*} \rho_{*} \cl \rightarrow \cl$ is surjective.

To see that $\rho_{*} \cl$ is torsion-free, note that since it is flat
and commutes with base change, it suffices to check the case where $R$
is a field and $ p = q = \pi = 0.$ In this case $\co_U = R[x,y]/sy$
and $\co_v = R[t,x]/tx$ and global sections of $\co(1)$ are $(g,f) \in
\co_U \oplus \co_V$ such that $g =sf.$ Moreover, $x(g,f)=(0,xf)$ and
$y(g,f) = (yg,0),$ hence if the ideal $(x,y)$ annihilates $(g,f),$
then $x(g,f) = y(g,f) = 0,$ and $xf=yg=0.$ Thus $f \in (t), g
\in (s),$ and this contradicts the fact that $g =sf$;
therefore, $ \rho_{*}
\co(1)$ has no associated primes of height one and is torsion-free.
\end{pf}

\begin{lemma}[Cornalba \cite{corn:theta}]
If $\cl$ is a line bundle on a semi-stable curve $f:X\rightarrow T$
such that $\cl |_{E} \cong \co _{E}$ for some exceptional component
$E$ of a special fibre $X_0,$ then there is an \'etale neighborhood
$T'$ of $0,$ such that $\cl|_T'$ on $X_{T'}$ is trivial in a
neighborhood of $E$ in $X_{T'}.$
\end{lemma}
\begin{pf}
There is an \'etale neighborhood $T''$ of $0$ in $T$ such that for
each irreducible component of the special fibre $X_{0}$ except $E,$
there is a section of $X''/T''$ that does not intersect $E$ but which
passes through the irreducible component (c.f. \cite[17.16.3]{ega4}).
Let $D$ be the divisor in $X$ corresponding to the sum of all these
sections, and note that $\cl (mD) |_{X_{0}}$ is generated by global
sections and has first cohomology group zero if $m$ is sufficiently
large.

So the natural map $ R^{1}f_{*}(\cl(mD)) \tensor k
\rightarrow R^1 f_* (X_{0},
\cl(mD) \tensor k) = 0$ is surjective, hence is an isomorphism, hence
$R^{1} f_{*} (\cl(mD))$ is zero on an open set $T'$ about $0$ in $T''$
(\cite[III.12.11]{H}).  This implies that on $T'$ the map $
R^{0}f_{*}(\cl(mD))
\tensor k \rightarrow
\Gamma (X_{0}, \cl (mD) \tensor k)$ is also surjective.
And $\cl(mD)$ is generated by global sections.  Thus $ \cl(mD)$ is
trivial on a neighborhood of $E,$ and on a sufficiently small
neighborhood of $E,$ $\cl(mD) \cong \cl.$
\end{pf}

As an immediate consequence of Cornalba's lemma we have the following
corollary.
\begin{corollary}
Proposition~\ref{tilden} holds for any line bundle $\cl$ which has
degree $1$ on the exceptional curve of the special fibre of
$\tilde{X}(p,q);$ namely,
\begin{enumerate}
\item $\rho_* \cl$ is flat over $R,$ it commutes with base change, and
$R^1 \rho_* \cl =0.$
\item $\Gamma(\tilde{X}, \cl) \cong \tilde{A}_1,$ and $\rho^* \rho_*
\cl \rightarrow \cl$ is surjective.
\item $\rho_* \cl$ is torsion-free of rank one.
\end{enumerate}
\end{corollary}

\subsection{Induced Maps}\label{indmap}

If $p$ and $q$ have the additional relations that $p^u = wq^v$ with
$u+v =r$ and $w
\in R^{\cross},$ then there is a canonical map from
$ \co_{\tilde{X}}(r)$ to $
\co_{\tilde{X}}.$  Namely, on $U$ it is $ \co(r) = (A[s]/ (ps-x, sy-q))
\cdot  {P}^{\tensor r} $ maps to
$A[s] / (ps-x, sy-q)$ via ${P}^{\tensor r} \mapsto wy^{v}.$ And on $V$
it is $ \co (r) = ( A[t]/ (p-xt, y-qt))
\cdot (\Xi)^{\tensor r}$ maps to $A [t] / (p-xt, y-qt)$ via
${\Xi}^{\tensor r}= s^{r} \cdot {P}^{\tensor r} \mapsto s^{r} wy^{v}
= x^{u}.$

The canonical map $ \rho^{*} \rho_* \co(1) \rightarrow \co (1) $
induces a map $ ( \rho^{*} \rho_* \co(1))^{\tensor r} =
\rho^{*}(\rho_{*} \co(1) ^{\tensor r})
\rightarrow \co(1)^{\tensor r} = \co (r),$ and the canonical
map $\co (r) \rightarrow \co_{\tilde{X}}$ gives a canonical map $
\rho^{*} (\rho_* \co(1)) ^{\tensor r} \rightarrow \co_{\tilde{X}}.$
And this induces a map on the push-forward by adjointness $$(\rho_*
\co(1))^{\tensor r} \rightarrow \rho_* \co_{\tilde{X}}= \
\co_{X}.$$
When we define various generalizations of a spin structure, the
best-behaved ones will be those which are locally isomorphic to those
induced from the canonical map $\co_{\tilde{X}}(r) @>>>
\co_{\tilde{X}}.$

\subsection{Local Structure of Torsion-Free Sheaves}\label{falt}

As in the previous section, we work with a stable curve $X/B,$ where
the base $B$ is the spectrum of a complete local Noetherian ring $R$;
the completion of the ring $\co_{X,\fp}$ at a singular point $\fp$ is
isomorphic to $A:=\rxy,$ and $\pi$ is an element of the maximal ideal
$\maxid$ of $R.$ A torsion-free sheaf $\ce$ corresponds to an $R$-flat
$A$-module, $E.$

Locally on $X,$ we can express a torsion-free sheaf obtained by
the contraction $\tilde X (p,q) @>>> X$ ($\ce$ is the
direct image of a line bundle $\cl$ of degree one on the exceptional
curve) in the following way.  $$E \cong
\Gamma(\spec A, \pi_{*} \co(1)) =
\tilde{A}_1=(A[\Xi,P]/(\Xi p-Px,  \Xi y-Pq))_{1} =
\{f \Xi+gP|f,g \in A\}.$$ It is easy to see that we can assume $f$ is
in $R[[x]],$ and $g$ is in $R[[y]].$ And so the map $\rho_{*}\co (1)
\rightarrow A^{\oplus 2},$ given by $(f\Xi+gP) \mapsto \left(
\begin{array}{c} fx+gp \\ fq+gy \end{array} \right)$ is a
well-defined homomorphism of $A$-modules.  Even if $p$ and $q$ are
zero divisors, if $f$ is in $R[[x]]$ and $g$ is in $R[[y]],$ then
$fx+gp=0$ implies that $f=0.$  Similarly, $fq+yg=0$ implies that
$g=0.$ So the map is injective.  We can, therefore, express $\rho_{*}
\co(1) $ as the image of the $A$-homomorphism $$\alpha (p,q)= \left (
\begin{array}{cc} x & p \\ q & y \end{array} \right) : A^{\oplus 2}
@>>> A^{\oplus 2}.$$

A result of Faltings shows that every rank-one torsion-free sheaf is
of this form.  Namely, let $E(p,q)$ be the image of $\alpha (p,q):
A^{\oplus 2} @>>> A^{\oplus 2},$ where $\alpha$ is the two-by-two
matrix $\left ( \begin{array} {cc} x & p \\ q & y \end{array}
\right),$ and $p$ and $q$ are, as before, elements of $R$ such that
$pq=\pi.$ We saw above that $E(p,q)$ is $R$-flat, and torsion-free.
When $p$ and $q$ are in $\maxid$ then $E(p,q)$ is a deformation of the
normalization of $A/\maxid A,$ i.e. of the unique (up to isomorphism)
non-free torsion-free sheaf $k[[x]] \oplus k[[y]]$ over the ring
$A/\maxid A \cong k[[x,y]]/xy.$ Faltings' result is the following.
\begin{theorem}[Faltings \cite{falt:torfree}]
Any reflexive $E$ of rank $1$ is isomorphic to an $E(p,q),$
for $p, q \in R$ with $ pq= \pi.$
\end{theorem}

The fact that
$\alpha \left( \begin{array}{c} y \\ 0 \end{array} \right) =
\alpha \left( \begin{array}{c} 0 \\ q\end{array} \right)$ and
$\alpha \left( \begin{array}{c} 0 \\ x \end{array} \right)=
\alpha \left( \begin{array}{c} p \\ 0 \end{array} \right)$
implies that for any $\alpha \left(
\begin{array}{c} f\\g \end{array} \right) \in E(p,q),$ we can assume
$f$ is in $R[[x]]$ and $g$ is in $R[[y]],$ and $E(p,q)$ is
$R$-isomorphic to $R[[x]] \oplus R[[y]]$ via the obvious
identification.  Homomorphisms and isomorphisms of $E(p,q)$'s can be
described by their lifts to $A^{\oplus 2}.$ Namely, any morphism of
$A$-modules $E \rightarrow F,$ with $E$ relatively torsion-free, can
be lifted to a morphism from $A^{\oplus 2}$ to $F.$ And a homomorphism
from $E$ to $E'$ with $E$ and $E'$ both torsion-free lifts to an
endomorphism of $A^{\oplus 2}.$ More exactly, the following holds.
\begin{proposition}[Faltings \cite{falt:torfree}]
If $ p \equiv q \equiv 0 \mod \maxid,$ then $ E(p,q)$ is isomorphic to
$ E(p',q')$ if and only if there exist $u,v \in R^{\cross}$ such that
$$ p'=upv^{-1}, \splice q' = vqu^{-1}.$$ In this case the isomporhism
is induced by the ``constant'' map $\left ( \begin{array} {cc} u & 0
\\ 0 & v \end{array} \right): A^{\oplus 2} @>>> A^{\oplus 2}.$
Moreover, writing a homomorphism $\Phi$ in $\Hom_{A} (E(p,q),
E(p',q'))$ as a lift to $A^{\oplus 2} @>>>A^{\oplus 2}$ given by
$\left(
\begin{array}{cc} \varphi_{+} & \psi_{+} \\ \psi_{-} & \varphi_{-}
\end{array}\right),$   $\varphi_+$ can be
taken to be in $R[[x]],$ $\varphi_-$ can be taken to be in $R[[y]],$
and the elements $\varphi_{+} (x)$ and $ \varphi_{-} (y)$ completely
determine $\psi_{+}$ and $\psi_{-}$ by the relations $$\psi_+ =
(\frac{p}{x}) (\varphi_+(x) - \varphi_+(0)) \splice \psi_- =
(\frac{q}{y}) (\varphi_-(y) -
\varphi_- (0)),$$ and are subject to the condition that $$ p'
\varphi_{-} (0) =
\varphi_{+} (0) p, \splice q' \varphi_{+} (0) = \varphi_{-} (0) q.$$

\end{proposition}

These results also hold for Henselian rings.  Suppose now that $R$ is
the Henselisation of a local ring of finite type over a field or an
excellent Dedekind domain, $\maxid$ is the maximal ideal of $R,$ $ \pi
\in \maxid,$ and $A$ is the Henselisation of $R[x,y]/(x y-\pi)$ at
$\maxid + (x,y).$ As before, for each pair $p, q
\in R$ with $ p q = \pi,$ define $E(p,q),$ and the
theorem is

\begin{theorem}[Faltings]

\begin{enumerate}
\item Any torsion-free $E$ of rank one over $A$ is
isomorphic to $E(p,q)$ for $p,q \in R$ and $p q
= \pi.$

\item If $p,q$ are in $\maxid,$  then $E(p,q)$ and
$E(p',q') $ are isomorphic if and only if there exist $u,v \in
R^{\cross}$ with $p'= upv^{-1}, q' = vqu^{-1}.$

\item  Suppose $I \subseteq R$ is a nilpotent ideal.  Modulo constant
automorphisms (given by $\left(\begin{array}{cc} u & 0\\ 0 & v
\end{array}\right)\in \text{GL}_2(R)$ with $up = pv, vq = qu$)
any automorphism of $E(p,q)/IE(p,q)$ lifts to an automorphism of
$E(p,q).$
\end{enumerate}
\end{theorem}

\section{Quasi-Spin Surves:  Local Structure}

A spin structure on a stable curve should have the property that where
the torsion-free sheaf is free it is isomorphic to the canonical line
bundle $\omega.$ A natural object to study, therefore, is a triple
$(X,\ce,b),$ with $X$ a stable curve, $\ce$ a relatively torsion-free
sheaf of rank one and degree $2g-2/r,$ and $b$ is a homomorphism of
$\co_{X}$-modules $$b:\ce^{\tensor r} \rightarrow
\omega_{\cx/T},$$ which is an isomorphism on the open set where $\ce$
is locally free.  Note that for smooth curves, such a triple is just
an $\rth$ root of the canonical bundle with an explicit isomorphism
of the \rth power of the bundle to the canonical bundle.  Later we
will need a few more conditions on these triples to get the ``right''
generalization of a spin curve, but we begin with these alone.

\subsection{A-Linear Homomorphisms of Tensor Powers}

To better understand these triples we need to study $A$-linear maps
$b: E^{\tensor r} @>>> A$ of the \rth\ tensor power of a rank-one
torsion-free $A$-module $E.$ As before, $A$ is an \'etale neighborhood
of the closed point defined by $(x,y) + \maxid A$ in
$\spec{\rxy}$ over the base ring $R,$ where $\maxid$
is a maximal ideal of $R$ containing $\pi, p,$ and $q.$

Any map $ b:E^{\tensor r} \rightarrow A$ that is $A$-linear, lifts to
a map $\tilde{b},$ thus
$$\begin{CD} A ^{\oplus 2^{r}} @>\tilde{b}>> A \\
@VV\alpha^{\tensor r}V @| \\ E^{\tensor r} @>b>> A
\end{CD}$$

Over $A[1/x]$ and over $A[1/y]$ the module $E$ is locally free, thus
over these rings any homomorphism $b:E^{\tensor r}
\rightarrow A$ will factor through $\text{Sym}^r(E),$ and its lift to
$(A^{\oplus 2})^{\tensor r}$ will factor through $\text{Sym}^r
(A^{\oplus r}).$  And since $A$ has no $(x,y)$-torsion, this holds in
general.  So if $f$ and $g \in A^{\tensor 2}$ are defined as $\left (
\begin{array} {c} 1 \\ 0 \end{array} \right)$ and $\left (
\begin{array} {c} 0 \\ 1 \end{array} \right)$ respectively, we only need to
describe $b_i :=b(f^{r-i} \tensor g^i)$ for each $ 0 \leq i \leq r$ in
order to completely describe $b$ and $\tilde{b}.$ We will, therefore,
denote $\tilde{b} : A^{2^{r}} \rightarrow E^{\tensor r}
\rightarrow A$ by the vector $(b_0, b_1, \dots, b_r).$

Now, since $\alpha \left( \begin{array}{c} p\\0 \end{array} \right) =
\alpha \left( \begin{array}{c} 0\\x \end{array} \right)$ and $\alpha
\left( \begin{array}{c} 0\\q \end{array} \right)  = \alpha \left(
\begin{array}{c} y\\0 \end{array} \right) ,$ we must have for all $i,$
$0 \leq i \leq r-1$
\begin{gather} \label{basicrelation}
pb_i =xb_{i+1}  \splice  yb_i =qb_{i+1}.
\end{gather}

Now the fact that the map must be an isomorphism off of the singular
locus of the underlying curve means that $b: E(p,q)^{\tensor r} @>>>
A$ must be an isomorphism on $A[1/x]$ and $A[1/y].$ Over $A[1/x]$ the
map $\tilde{b} = (b_0,b_1,\dots,b_r)$ is completely determined by
$b_0,$ namely for any $i$ we have $x^ib_i = p^ib_0,$ and thus $$b_i =
b_0 p^i/x^i.$$ For $\tilde b$ to be surjective over $A[1/x]$ we must
have that $b_0$ is invertible in $A[1/x].$ Similarly, $b_r$ is
invertible in $A[1/y].$

On the special fibre, since $p$ and $q$ are in $\maxid,$ we
have $\tilde b \equiv (\bar b_0,0,0,\dots,0,\bar b_r)
\pmod{\maxid}.$  Here $\bar z$ denotes the image of $z$ modulo $\maxid.$
And $\bar b_0=x^u\bar\beta_0$ and $\bar b_r=y^v\bar\beta_r$ for some
$\bar\beta_0$ invertible in $(A/\maxid A)[1/x],$ but not in the ideal
$(x),$ and for some $\bar\beta_r$ invertible in $(A/\maxid A)[1/y],$
but not in $(y).$ This makes the length of the cokernel of $\bar b$
equal to $u+v-1.$ If $u$ (or similarly $v$) is zero, then over the
ring $A[1/y]$ the fact that $b_0 = b_r q^r/y^r \equiv 0
\pmod{\maxid}$  implies that $\bar\beta_0 = \bar b_0 \equiv 0
\pmod{\ann_{A/\maxid A}(y))},$ i.e. modulo $(x).$  But this
implies that $\bar\beta_0 = 0,$ and that is a contradiction.  Hence we
have proven the following.
\begin{proposition} \label{b}
Any pair $(E,b)$ which is not free is of the form $(E(p,q),b)$ with
$p$ and $q$ in $\maxid,$ and $b$ lifts to $\tilde b=
(b_0,b_1,\dots,b_r),$ where, modulo the ideal $\maxid A$ we have $\bar
b_0=x^u\bar\beta_0$ and $\bar b_r=y^v\bar\beta_r$ for some
$\bar\beta_0$ invertible in $(A/\maxid A)[1/x],$ but not in the ideal
$(x),$ and for some $\bar\beta_r$ invertible in $(A/\maxid A)[1/y],$
but not in $(y).$ Moreover, $u$ and $v$ must both be at least one.
\end{proposition}

The final condition we will need on triples $(X, \ce, b)$ to generalize
spin curves is the condition that the two constants $u$ and $v$ in
Proposition~\ref{b} must sum to $r.$ In other words, the length of the
cokernel of $b$ at each singular point of each fibre must be $r-1$ if
the sheaf $\ce$ is not free there.  Without this condition there would
be too many possible triples for the associated stack to be separated.
We now have the first generalization of spin structures on stable
curves.
\begin{defn}
A {\em quasi-spin structure} on a stable curve $X/T$ is a pair $(\ce,b)$
where $\ce$ is relatively torsion-free of rank one, and $b$ is a
homomorphism of $\co_{X}$-modules $$b:\ce^{\tensor r} \rightarrow
\omega_{X/T}$$ to the canonical dualizing sheaf, such that \begin{enumerate}
\item $\ce$ has degree $(2g-2)/r.$
\item $b$ is an isomorphism on the open set where $\ce$ is not free.
\item For each closed point $t$ of the base $T,$ and for each singular point
$\fp$ of the fibre $\cx_{t}$ where $\ce$ is not free, the length of
the cokernel of $b$ at $\fp$ is $r-1.$
\end{enumerate}
\end{defn}
\begin{defn}
A {\em quasi-spin curve} is simply a stable curve with a quasi-spin
structure.
\end{defn}

In the special case that $r=2,$ the requirement that the cokernel of
the spin structure map be supported on the singular locus of $\ce$ is
enough to guarantee that the length of the cokernel is at least one at
all singular points.  The condition on the total degree of $\ce$ can
be seen to guarantee that the length of the cokernel is at most (and
hence exactly) one at all singular points.  Moreover, these conditions
are equivalent to the condition that the map $b$ induce an isomorphism
$\ce \irightarrow
\Hom_{\co_{\cx}}(\ce,\omega_{\cx})=\ce^{\vee}\tensor \omega_{\cx} .$

\subsection{Power Series Expansions}

Any $A$-linear map $ b=(b_{0}, \dots, b_{r})$ as above, over the complete
local ring $A=\hat{\co}_{\cx,\fp}$ has a power series expansion $b_{i}
= \sum_{n\geq 0} b_{in} x^{n} +
\sum_{m>0} b_{i,-m} y^{m}.$  And the relations $p^i b_0 = x^i b_i$ and
$q^{r-i} b_r = y^{r-i} b_i$ imply that $$p^i b_{0,n+i} =b_{i,n}$$ for
$n \geq 0,$ and $$ b_{i,-m} = q^{r-i} b_{r, -m-(r-i)}$$ for $ m \geq
0.$ And in particular $$ p^{j} b_{0,j} = q^{r-j} b_{r, j-r}$$ for all
$ j,0 \leq j \leq r.$

Moreover, if $b$ induces a quasi-spin structure on the central fibre,
then there are $ u,v \in \Bbb Z^{+}$ such that  $\pmod \maxid$
$$\bar{b}_{0} = x^{u} \bar{\beta}_{0} \splice
\bar{\beta}_{0} \in (\bar{A}_{x})^{\cross}.$$
This implies that $\bar{b}_{0} = \sum_{n \geq u} \bar{b}_{0,n} x^{n}$
with $
\bar{b}_{0,u} \neq 0,$ hence $b_{0,u} $ is not in $
\maxid$ and is invertible in $R.$
Similarly, $b_{r,-v} \in R^{\cross}.$ So, in particular, $ p^{u}=q^{v}
b_{r,-v}/b_{0,u}.$ Letting $ w = b_{r,-v}/b_{0,u} \in A^{\cross},$ we
have the relation $$ p^{u} = q^{v}w.$$

In the special case that $\pi$ is not a zero divisor, the relations $
p^{i} b_{0,i} = q^{r-i} b_{r, i-r}$ for $0 \leq i \leq u$ imply that
$$ b_{0,i} = w^{-1} p^{u-i} q^{u-i} b_{r,i-r} =
\frac{\pi^{u-i}}{w} b_{r, i-r}. $$ Similarly, $$ b_{r, i-r} = w
\pi^{i-u} b_{0, i}, \text{ for } u \leq i
\leq r.$$ But even when $\pi$ is a zero divisor $$b_{0,0} = \pi^u b_{r,
-r} \commasplice b_{r,0} = \pi^v b_{0,r} \commasplice b_{0,u} = w
b_{r, -v},$$ and $$b_{0,i}= \frac{\pi^{u-i}}{w} b_{r, i-r} + \sigma_i
\text{ for } 0 < i < u, \splice b_{r, i-r} = w \pi^{i-u}b_{0,i} + \sigma_i
\text{ for } u<i < r.$$

The ``bad'' terms $\sigma_i$ are all nilpotent elements. On the one
hand, for any prime ideal $\fp \in \spec R$ such that $p$ (and hence
$q$) is in $\fp,$ we have that $ b_i \equiv 0 \mod
\fp$ for $ 0 <i<r.$ And $b_0 \equiv x^u \beta,$ and $b_r \equiv y^v
\gamma,$ with $\beta$ and $\gamma$ invertible elements of $R[[x]]$ and
$R[[y]]$ respectively.  Accordingly, $$b_{0,i} \in \fp \text{ for } 0
< i <u, \splice b_{r, i-r} \in \fp \text{ for } u <i < r,$$ and thus
$\sigma_i \in \fp$ for $0 <i < r$ ($\sigma_u$ is obviously zero).  On
the other hand, if $p$ (and therefore $q$) is not in $\fp,$ then $p$
and $q$ are not zero divisors in $R_{\fp}$ and hence, as demonstrated
before, $\sigma_i \in \fp$ for $ 0 <i <r .$ Thus $\sigma_i$ is
contained in the nilradical of $R$ for every $i.$ And in particular,
for any quasi-spin structure over a reduced, complete local ring the
relations
\begin{gather}\label{ssrelations}
 b_{0,i} = \frac{\pi^{u-i}}{w} b_{r, i-r} \text{ for } 0 \leq i \leq
u \\ b_{r, i-r} = w \pi^{i-u} b_{0,i} \text{ for } r \geq i \geq
u \notag
\end{gather}
hold.

\begin{proposition}
When the relations (\ref{ssrelations}) hold, we can write $b_{0}$ and
$b_{r}$ as the following products: $$b_0=ax^u \splice b_r = awy^v
\text{ for some } a \in A^{\cross}.$$  In particular, given $u,v,$ and
$w,$ the fact that the specified relations (\ref{ssrelations}) hold
means that $b$ is completely determined by $a \in A^{\cross}.$
\end{proposition}
\begin{pf}
$$ b_{0} = x^{u} \sum_{n \geq 0} b_{0, n+u} x^{n} +
1/w \sum_{u>m\geq 0} \pi^{u-m} b_{r, m-r} x^{m} + q^{r} \sum_{l>0}
b_{r,l-r} y^{l}$$ And thus $$ a = \sum_{n \geq 0} b_{0, n+u}
x^{n} + 1/w \sum_{m > 0 } b_{r, -m-v} y^{m}.$$ The calculation is
similar for $b_{r}.$
\end{pf}

\begin{proposition}
If the relations (\ref{ssrelations}) hold on $\tilde b,$  then
$b$ is actually the map induced on $E(p,q)^{\tensor r} = \left (
\rho_{*} \co_{\tilde{X}(p,q)}(1) \right )^{\tensor r}$ as in
Section~\ref{indmap}.  Moreover, the relations (\ref{ssrelations})
hold for (the lift to $A^{\oplus 2^r}$ of) the map induced from
$\co_{\tilde{X}(p,q)}(1)$ for any $p$ and $q$ in $\maxid_R$ with
$pq=\pi.$
\end{proposition}

\begin{pf}
To see this out explicitly the map $$ (A^{2})^{\tensor r}
\tilde{@>>>} A^{2^{r}} @> \alpha^{r}>>(E(p,q))^{r}
\tilde{@>>\psi^{r}>} (\rho_{*} \co(1))^{r} @>>\varphi> \rho_{*} \co(r) @>>
\gamma> \co_{X} = A.$$
The first map is just $$ \left (
\begin{array}{c} f_{1} \\ g_{1} \end{array} \right )\tensor \dots
\left ( \begin{array}{c} f_{r} \\ g_{r} \end{array} \right )
\mapsto \alpha \left ( \begin{array}{c} f_{1} \\ g_{1}
\end{array} \right ) \tensor \dots \tensor \alpha \left ( \begin{array}{c}
f_{r} \\ g_{r} \end{array} \right ).$$ The map $\psi$ is given by $
\psi : \alpha \left (
\begin{array}{c} f \\ g \end{array} \right ) \mapsto (sf+g, f+tg),$ so
$$\psi^{\tensor r} : \alpha \left ( \begin{array}{c} f_{1} \\ g_{1}
\end{array} \right ) \tensor \dots \tensor \alpha \left (
\begin{array}{c} f_{r} \\ g_{r} \end{array} \right ) \mapsto (sf_{1}
+ g_{1}, f_{1} + tg_{1}) \tensor \dots \tensor (sf_{r} + g_{r}, f_{r} +
tg_{r}).$$
The map $\varphi$ is given by
$ \varphi : (h_{1}, k_{1}) \tensor \dots \tensor (h_{r},
k_{r} ) \mapsto ( h_{1} h_{2} \dots h_{r}, k_{1} k_{2} \dots
k_{r})$
and $\gamma$ is $ \gamma: (h, k) \mapsto w y^{v} h = x^{u}k \in A.$

So the composite map is $$\left ( \begin{array}{c} f_{1} \\ g_{1}
\end{array} \right ) \tensor \dots \left ( \begin{array}{c} f_{r} \\
g_{r} \end{array} \right ) \mapsto x^{u} \prod_{1 \leq i \leq r}
(f_{i} + tg_{i}),$$ but this is just the map $$ b = (x^{u}, tx^{u}
\dots, t^{r} x^{u}) = ( x^{u}, px^{u-1}, \dots, p^{u}, p^{u}t, \dots,
p^{u} t^{v}). $$ And $ p^{u} = wq^{v},$ and $ qt=y,$ so $$ b = (x^{u},
px^{u-1}, \dots, p^{u}, wq^{v-1}y, \dots, w y^{v}).$$

This composite map depends only on the choice of isomorphism $
\Gamma(\tilde{X}, \co_{\tilde{X}}) \irightarrow A,$ and any
element $a$ in $A^{\cross}$ induces an automorphism $A,$ so any
quasi-spin structure with the additional relations (\ref{ssrelations})
is actually an induced map.
\end{pf}

In the special case when $r=2$ every singularity has $u=1,$ and thus
$\sigma_0 =\sigma_1=\sigma_2=0.$ Therefore all quasi-spin structures
are actually locally isomorphic to an induced structure.

We can always map $ E(p,q) \tilde{\rightarrow} E(p',q')$ with $p' =
\lambda p, q' = \lambda^{-1} q.$  So if $ p^{u} = wq^{v},$ then
${p'}^{u} = \lambda^{u} p^{u} = \lambda^{u} wq^{v} = \lambda^{r}
w{q'}^{v}.$ So $ w \rightarrow \lambda^{r} w.$ And $w$ has an
\rth root in $k$ if the base ring
$R$ has its residue field $k$ algebraically closed, hence by the
step-by-step method the \rth root of $w$ will lift to all of $R,$ i.e.
if $\maxid I = 0$ there is an $i \in I$ such that $ w + i =
\lambda^{r},$ which implies that $(\lambda - \frac {i}{r
\lambda^{r-1}})^{r} = w.$ Thus we may assume that if the central fibre
has residue field $k =
\bar{k},$ then $w$ can be taken to be one, and, in general, $R$ has an
\'etale cover on which we can take $w$ to be one.  Accordingly, all
$b$'s for which the relations (\ref{ssrelations}) hold are determined,
\'etale locally, by $ p,q,u,v,$ and an element of $A^{\cross}.$

\subsection{Behavior of the Cokernel Under Deformation}
\label{sec:coker}
Because of the condition on the cokernel of the quasi-spin map
$b,$ we need to understand the way that the length of the cokernel
changes under deformation.

\begin{defn}
An $\co_{X}$-linear map $b:\ce^{\tensor r} @>>> \omega$ from the \rth\
tensor power of a rank-one torsion-free sheaf $\ce$ to the canonical
bundle of a curve $X$ over $k$ is said to {\em have good cokernel} if
\begin{enumerate}
\item the cokernel is supported on the singularities of $X,$ and
\item for each point $\fp$ of the support of the cokernel $C$ of $b$
$$ \len_{\fp}C = r-1 .$$
\end{enumerate}
\end{defn}

\begin{lemma}
If the cokernel of $b$ is supported on the singular locus of $\ce,$
then the property of having good cokernel is stable under
generization.
\end{lemma}
\begin{pf}
It is enough to consider the case where $R$ is a complete local ring,
$E \cong E(p,q)$ is an $A$-module, with $A=\rxy,$ and
$\tilde{b}=(b_0,\dots,b_r):A^{2^r} @>>>A$ is a lifting of the map
$b:E^{\tensor r} @>>> A.$ We can assume that on the special fibre
$\overline{b}_0 \in A/\maxid A$ equal to $x^{i}\overline{\beta}_0,$
with $\overline{\beta}_{0}$ an invertible element of $( A/\maxid).$
Similarly, $\overline{b}_r \in A/\maxid A$ with $\overline{b}_r =
y^{j} \overline{\beta}_r,$ and $\beta_{r}$ invertible in $A.$

 Now for any map $R @>>> K$ of $R$ into a field, we have the following
possible cases.
\begin{enumerate}
 \item $\pi$ does not map to zero in $K.$ In this case, the cokernel
is actually zero because \spec{A\tensor K} is regular.  \item $\pi$
maps to zero, but at least one of $p$ and $q$ does not.  In this case
again the cokernel of $b$ is zero.  \item $\pi$ and $p$ and $q$ all
map to zero.  This is the only interesting case.  We have $A\tensor K
\cong K[[x,y]]/xy$ and $\tilde{b}_K = (b_0,0,\dots,0,b_r).$ Now
$$b_{0} = x^i \beta_0 + d_0 \splice b_{r} = y^j \beta_r + d_r$$ with
$d_{0}=x^{l}\ce_{0}$ and $d_{r}=y^{m} \ce _{r},$ such that $\ce _{0}$
is in $\maxid A$ but not in $(x),$ and $\ce _{r}$ is in $\maxid A$ but
not in $(y).$ If, on the one hand, $l$ is larger than $i-1,$ then
$$x^{i}=b_0/(\beta_0 + x^{l-i}\ce_{0}).$$ The term in the denominator
is invertible because $\beta_0$ is invertible, and $d_0$ is in the
maximal ideal $\maxid A.$ If, on the other hand, $l$ is less than $i,$
then $$x^i = \frac{-x^l \ce_0 }{\beta_0},$$ and similarly for
$y^j.$  In either case  $$K[[x,y]]/(xy,x^i,y^j) \cong
K < 1, x, x^2, \dots, x^{i-1}, y, y^2, \dots, y^{j-1}>$$ surjects onto
$$K[[x,y]]/(xy,b_0,b_r)=A \tensor_R K/\im{b}.$$
\end{enumerate}
So the length of the cokernel will be either zero (cases 1 and 2) or
bounded above by $ i+j-1=r-1$ (case 3).

Thus the length of the cokernel can only decrease under generization,
but the degree of $\ce_K$ on $\cx_K$ must be $(2g-2)/r = \deg
\theta^{\natural} \ce + \delta,$ where $\delta$ is the number of
singularities of $\cx_K,$ and $\theta: \cx^{\nu}_K \rightarrow \cx_K$
is the normalization of $\cx_K$ at the singularities of $\ce_K.$  On
the other hand, since the cokernel of $b$ is supported on the singular
set of $\ce,$ we have that $\theta^{\natural} b$ factors
$$\theta^{\natural} \ce^{\tensor r}_K
\irightarrow \theta^{*}
\omega_{\cx_K} (-\sum u_{\fp}\fp^+ -\sum v_{\fp}\fp^-) \hookrightarrow
\theta^{*} \omega_{\cx_K},$$ where the sum is taken over all $\fp$ in the
singular set of $\ce_K,$ $\theta^{-1}(\fp) = \{\fp^{\pm}\},$ and for
each $\fp,$ $u_{\fp}+v_{\fp} -1=\len_{\fp} (\coker(b)) \leq r-1.$ So
$$\deg \ce_K = (2g-2)/r = \bigg(2g-2-\sum_{\fp} (u_{\fp} +
v_{\fp})\bigg)/r +
\delta,$$ which will be strictly greater than $(2g-2)/r$ unless at each
singularity of $\ce_K$ the cokernel of $b$ has length $r-1.$ Thus the
property of having good cokernel is stable under generization.
\end{pf}

\begin{proposition}
\label{opencoker}
Given $ b: \ce^{\tensor r} \rightarrow \omega$ on $f:\cx @>>>T$ (with
the cokernel of $b$ supported on the discriminant locus) the set of $t
\in T$ such that $b_t$ has good cokernel is open in $T.$
In other
words, the functor of $T$-schemes $$F_{b} (T')= \begin{cases} \{1\} &
\text{if $b$ has good cokernel at every geometric point of $T'$} \\
\emptyset & \text{if there exists $\overline{t} \in T'$ where
$b$ does not have good cokernel}
\end{cases}$$
is an open subfunctor of the trivial functor $T' \mapsto \{1\}.$
\end{proposition}

\begin{pf}
It suffices to show the complement of the set is closed.  And since
the previous lemma shows this complement is stable under
specialization, it suffices to show that the complement is
constructible.  Let $P_{m}$ be the property of a geometric point
\tbar\ of $T$ that $C:= \coker (b)$ has a point of its support over
\tbar\ where $C$ has length $m.$ The set we want to show is
constructible is the set $\frak{T}_m :=\{ t \in T |
\tbar
\text{ has } P_{m}\}.$  Actually, the set we are really looking
for is $$\bigcup\begin{Sb}{0 < m < r-1}\\ r - 1 < m < N \end{Sb}
\frak{T}_{m},$$ for some very large $N.$  $N$ can be taken to be
finite because the degree of the sheaves is fixed, and the sum over
all points in a given fibre of the length of the cokernel is bounded,
and this bound is determined by the number of singular points and the
degree.  Moreover, the number of singular points is bounded as a
function of the genus of the underlying curve, so this number $N$ can
be chosen independently of the specific family.

Now to show constructibility we only need to consider one $m$ and one
irreducible component of the discriminant locus, say $D_{0},$ of $\cx$
and its image $\rho (D_{0} ) = T_{0},$ i.e. we only need to show that
$\frak{T}_m$ is constructible in $T_0.$  And it is enough to assume
$T_{0}$ is reduced and irreducible.  Since $D_{0} $ is proper over
$T_{0},$ the semi-continuity theorem shows that $\frak{T}_m$ is
constructible for any $m \neq r-1.$
\end{pf}

\section{Local-to-Global Calculations}
\subsection{Log-Structures}

The well-behaved quasi-spin curves, i.e. those for which the relations
(\ref{ssrelations}) hold locally, also give a compactification of the
moduli space of spin curves, and their local moduli spaces are much
easier to describe than those of general quasi-spin curves.  But in
order to formalize the notion of ``well-behaved'' we need to choose
local coordinates for the whole curve in such a way that our
constructions make sense globally.

 From the deformation
theory of stable curves, we know that the complete local ring
$\hat{\co}_{\cx,x}$ over $\hat{\co}_{T,t}$ is of the form
$\hat{\co}_{\cx,x} \cong \hat{\co}_{T,t} [[x,y]]/(xy-\pi)$ for some $\pi
\in \hat{\co}_{T,t}.$  And on some \'etale neighborhood $T'$ of $t,$
 the induced curve $\cx \cross_T T'$ has the
following additional structure:  on an \'etale cover $\cx'$ of ${\cx}
\cross_T T',$ there are
sections $x$ and $y$ in $\co_{{\cx}'}$ such that
\begin{enumerate}
\item $xy=\pi \in \co_{T,t}.$
\item The ideal generated by $x$ and $y$ has the discriminant
locus of $\cx/T$ as its associated closed subscheme.
\item The obvious homomorphism $\big(\co_{T,t}
[x,y]/(xy-\pi)\big) \rightarrow
\co_{{\cx}',x}$ induces an isomorphism on the completions
$\big(\hat{\co}_{T,t} [[x,y]]/(xy-\pi)\big) \irightarrow
\hat{\co}_{{\cx}', x}.$
\end{enumerate}

Such a collection of data ($\cx',T',x,y,\pi$) is what we need locally.
But this data is not uniquely determined; it is only determined up to
the equivalence relation generated by the operations
\begin{enumerate}
\item pullback to \'etale covers.
\item change by units: namely $x' = \tilde{u}x, y'=\tilde{v}y, {\pi}'
= \tilde{w} \pi$ with $\tilde{u}, \tilde{v} \in \co^*_{{\cx}'},$ and
$\tilde{u} \tilde{v} = \tilde{w} \in \co^*_{T'}.$
\item switching branches: namely, interchanging $x$ and $y.$
\end{enumerate}

A log structure is a way of choosing these local data coherently.
\begin{defn}
 A {\em log structure} for $\cx/T$ is given by \'etale covers $\cx'$
and $T'$ of $\cx$ and $T$ $$\begin{CD} \cx @<<< {\cx}' \\
 @VVV @VVV \\ T @<<< T'
\end{CD}.$$ And for each irreducible component of the singular locus of
${\cx}'$ a choice of $\pi \in \co_{T'}$ and a choice of $x$ and $y$ in
$\co_{{\cx}'}$ with the three properties listed above, and with
descent data related to the equivalence relation.  Namely, on ${\cx}''
= {\cx}' \cross_{\cx} {\cx}'$ over $T'' = T' \cross_T T'$ with
projection maps $pr_1$ and $pr_2,$ there are $1$-cocycles $u,v$ in
$\co^*_{\cx''},$ and $w$ in $\co^*_{T''}$ such that: $\text{pr}_2^*
(x) = u\text{pr}_1^*(x), \text{pr}_2^* (y) = v\text{pr}_1^* (y),$ and $uv
=w,$ and $\text{pr}_2^*(\pi) = w\text{pr}_1^* (\pi)$ with the cocycle
condition that on ${\cx}''' = {\cx}' \cross_{\cx}{\cx}'
\cross_{\cx}{\cx}',$ $u,v,$ and $w$ are all compatible with the
different projections, i.e. $ \text{pr}^*_{12}(u)\text{pr}^*_{23}(u) =
\text{pr}^*_{13}(u)$ and so forth.
\end{defn}

As in the local case, we also impose the equivalence relation on the
log structures generated by pullback to \'etale covers and by change
by units compatible with the descent data; namely, two log structures
$(\cx', T', x,y,\pi)$ and $(\cx', T', x', y', \pi ')$ are equivalent
if there exist $\tilde{u}, \tilde{v}$ in $\co^*_{\cx'}$ and
$\tilde{w}$ in $\co^*_{T'}$ such that $x' = \tilde{u}x,
y'= \tilde{v} x, \pi ' = \tilde{w} \pi, $ with $\tilde{u} \tilde{v} =
\tilde{w}$ and if $(u,v,w)$ and $(u',v',w')$ are the cocycles
corresponding to the two log structures, the units $\tilde{u},
\tilde{v} \splice \tilde{w}$ must be compatible with them as
well, namely $u' = (\text{pr}^*_1 (\tilde{u})/\text{pr}^*_2
(\tilde{u})) u,$ and $ v' = (\text{pr}^*_1(\tilde{v})/\text{pr}^*_2
(\tilde{v}))v,$ and $ w' = (\text{pr}^*_1 (\tilde{w})/\text{pr}^*_2
(\tilde{w}))w.$

As discussed above, given any two log structures with distinguished
branches $(x)$ and $(y)$ we will have the relations $x' = \tilde{u} x$
and $y' = \tilde{v} y,$ etc.  And they will be almost equivalent,
namely $\text{pr}^*_1 (x)
(u'-(\text{pr}^*_1(\tilde{u})/\text{pr}^*_2(\tilde{u}))u)=0$ and so
forth; thus if $\pi$ is not a zero divisor (and hence $x$ and $y$
also) all log structures are equivalent.  In particular, since a
versal deformation of the curve has no zero divisors it has a unique
log structure.

Switching of branches $(x\mapsto y, y\mapsto x)$ and switching of
double points (i.e. interchange the different $\pi_i$) results in an
action of the $n\th$ symmetric group ($n$ is the number of
double points) and the group $(\Bbb Z/2\Bbb Z)^n$ on the log
structures.  But for our purposes this is not a problem, namely we are
interested in expressing $\ce$ as an $E(p,q)$ and this switching just
interchanges $p$ and $q$ or the different $\pi_i.$ So given a log
structure on a stable curve, we can use the methods of Faltings to
describe rank-one torsion-free sheaves, namely any such sheaf $\ce$ is
isomorphic to an $E(p,q),$ and the results on homomorphisms and
isomorphisms still hold.

In general the choice of a log structure is unique up to the
automorphisms $x \mapsto ux,$ $y\mapsto vy,$ and $\pi
\mapsto w\pi,$ for $uv=w,$ but locally this might not be all of the
automorphisms of the henselization of the ring $R[x,y]/(xy-\pi).$ In
other words, on a curve $C \rightarrow B$ we might have different log
structures induced by different maps of $B$ to the versal deformation.
Nevertheless, we can get around this by considering the problem
globally instead.  Namely let $\cs/\crr = \mgbar$ be a presentation of
the stack of of stable curves, i.e. $\cs$ is \'etale over \mgbar, and
$\crr$ is the \'etale equivalence relation (\Isom).  $\crr$ is smooth
and has no zero divisors, so the two pullbacks to $\crr$ of the
universal curve with its unique log structure over $\cs$ are
canonically isomorphic.  Hence any curve over any base has a canonical
log structure induced by the unique log structure on the universal
curve over $\cs.$

Note that given a choice of $p$ and $q$ in $\co_T,$ the descent data
for the canonical log structure determine gluing data for the various
blowings up.  Thus the techniques of Section~\ref{geom} yield a
globally-defined semi-stable curve $\tilde{\cx}(p,q)$ over $\cx,$ a
rank-one torsion-free sheaf $\ce(p,q) = P_{*} \co (1)$ on $\cx,$ and a
canonical map $b:\ce^{\tensor r} @>>> \cm,$ for some line bundle
$\cm.$

\subsection{Spin Curves and Pure-Spin Curves}

We can now define our ``good'' quasi-spin curves using the canonical
log structure.
\begin{defn}
 A {\em spin structure} on an arbitrary
stable curve $\cx/T$ is a pair $(\ce, b),$ where $\ce$ is
relatively torsion-free of rank one with degree $(2g-2)/r,$ and $b$ is
a morphism of $\co_{X}$-modules $$b:\ce^{\tensor r} \rightarrow
\omega_{\cx/T},$$ which is an isomorphism on the open set where $\ce$
is locally free, and such that via the canonical log structure on
$\cx/T,$ the sheaf $\ce$ is isomorphic to $E(p,q)$ for some $p$ and
$q$ in $\co_T$ with $pq=\pi,$ and the homomorphism $b:\ce^{\tensor r}
\rightarrow \omega$ is the canonical induced morphism.
\end{defn}

An even stronger condition that we can impose on the spin curves is
that $p$ and $q$ be such that $p=t^v$ and $q=t^u$ for some $t$ in
$\co_T.$ Spin curves that have this property will be called {\em pure
spin curves.}  Pure-spin curves also compactify the smooth spin curves,
and they have an especially well-behaved local structure, as we will
see later.

\section{Deformation Theory}\label{deftheory}

Given a spin structure $(\bar{\ce}, \bar{b})$ on a curve $\bar{\cx}$
over an Artin local ring $\bar{R}$ with residue field $k,$ and given a
deformation $R$ of $\bar{R},$ namely $\bar{R}= R/I$ with $I^2 =0,$ we
want to study deformations of $(\bar{\cx}, \bar{\ce}, \bar{b})$ to
spin curves and quasi-spin curves over $R.$

First, do this locally.  For $\bar{A}:= \bar{R}[[x,y]]/(xy-\bar{\pi})$
and $(\bar{E},\bar{b}) =(E(\bar{p}, \bar{q}), \bar{b}),$ a spin
structure on $\bar{A},$ we want to lift $\bar{A}$ and $(\bar{E},\bar{b}).$
But any lift corresponds to a choice of $\beta,$ $P$ and $Q$ such that
$P^u=Q^v$ and an isomorphism $({\overline{E ({P}, {Q})}},\bar{\beta})
\irightarrow (E(\bar{p},\bar{q}),\bar{b}).$  Here ${\overline{E
(P,Q)}}$ is the module $E(P,Q)/I\cdot E(P,Q) = E(\bar P,\bar Q)$
induced by pulling back $E(P,Q)$ along the canonical map $\spec R
\leftarrow \spec{\bar R},$ and the map $\bar{\beta}$ is the canonical map
induced from $\beta$ on $\overline{E(P,Q)}.$ By Faltings' theorem,
isomorphisms over $\overline{R}$ are of the form $\bar{\Phi} = \left
(\begin{array}{cc} \bar{\zeta} & 0 \\ 0 &
\bar{\xi} \end{array} \right) $ with $\bar{\zeta}, \bar{\xi} \in
\bar{R}^{\cross}.$ These lift to isomorphisms $\Phi = \left(
\begin{array}{cc} \zeta & 0 \\ 0 & \xi
\end{array} \right) $
 for lifts $\zeta,\xi \in R^{\cross}$ of $\bar{\zeta}$ and
$\bar{\xi},$ and such lifts always exist in $R.$ Thus any local spin
structure is given simply by a choice of $P$ and $Q$ in $R$ such that
$\bar P = \bar p$ and $\bar Q = \bar q$ and a choice of $\beta$ such
that the induced map $\bar{\beta}$ on $E(\bar p,\bar q)^{\tensor r} =
\bar E^{\tensor r}$ differs from $\bar b$ only by an automorphism of
$\bar E.$ In particular, $\bar{\beta} = \bar{a} \bar{b}$ with
$\bar{a}$ in $\bar{A}^{\cross},$ and thus $\beta$ is uniquely
determined by an element $a \in A^{\cross},$ i.e.
$\beta=a(x^u,px^{u-1},\dots, y^v).$ This describes the local
deformations.

Any combination of local lifts will patch together into a global one.
This is due to the fact that if $\bar{\sigma}$ is a section of
$\co^*_{\bar{\cx}}$ and $\gamma$ is a section of $\co^*_{\cx}$
inducing $\bar{\gamma}$ in $\co^*_{\bar{\cx}}$ such that
$\bar{\sigma}^r = \bar{\gamma},$ then $\bar{\sigma}$ lifts uniquely to
a section $\sigma$ of $\co^*_{\cx}$ such that $\sigma^r = \gamma.$
This is easy to check.  We can now use {\em fpqc} descent to lift
$(\bar{\cx}, \bar{\ce}, \bar{b}).$  Namely, a choice of $P$ and $Q$ for
each singularity of $\bar{\cx}$ still allows numerous choices of $\cx$
deforming $\bar{\cx}.$  And given such an $\cx,$ we need to construct
a pair $(\ce,b)$ extending $\bar{\ce}$ and $\bar{b}.$  On $U,$ the
complement of the discriminant locus, the line bundle $\bar{\ce}$
extends uniquely to a line bundle $\ce$ that is an \rth\ root of
$\omega.$  Given an extension $E(P,Q)$ at each singularity (i.e. at
$\spec{\hat{\co}_{\cx, x_i}}$), we have a covering datum induced
by the unique lift of the covering datum on $\bar{\cx}$ that makes
$\bar{\ce}$ an \rth\ root of $\bar{\omega}.$  This datum is actually a
descent datum because of the uniqueness of \rth\ root lifts.  And
since all {\em fpqc} descent data for coherent sheaves are effective,
we have the desired $(\ce,b)$ on $\cx$ extending $(\bar{\ce},\bar{b}).$

In fact the universal deformation of a spin curve $(X,E, \beta)$ over
a field $k$ is the obvious formal spin curve $(\tilde{\cx}, \ce, b).$
Here $\tilde{\cx}$ is the pullback of the universal deformation $\cx
\rightarrow
\frak{o}_{k} [[t_1,\dots, t_{n}]]$ of the curve $X/k$ along the
homomorphism $$\frak{o}_k [[t_1, \dots, t_n]]
\rightarrow \frak{o}_k [[P_1, Q_1, \dots P_l, Q_l, t_{l+1}, \dots,
t_{n} ]]/(P_i^{u_i}-Q_i^{v_i})$$ via $t_i \mapsto P_iQ_i$ for $i
\leq l.$  $u_i$ and $v_i$ are determined by the map $\beta$ at each
singularity of $E$ on the central fibre.

A particularly useful corollary of this is the following.
\begin{proposition}\label{smoothdef}
All quasi-spin structures over a field have a deformation to a smooth
spin structure.
\end{proposition}

\section{Isomorphisms}

\subsection{Isomorphisms of Spin Structures over a Field}

Since spin structures, quasi-spin structures and pure-spin structures
are all the same over a field, the study of isomorphisms over a field
is fairly simple.  Any two spin structures on $ X/k,$ say $(\ce,b)$
and $(\ce', b'),$ which are singular at the same points and are the
same on $X^{\nu}$ ($X^{\nu} @>{\theta}>> X$ is the normalization of
$X$ at the singularities of $\ce$) via $\theta^{\natural},$ must be
isomorphic on $X.$  For a given $\ce$ on $X,$ any two spin structure
maps $b$ and $b'$ are the same if and only if $\theta^{\natural} b=
\theta^{\natural} b',$ and this is true if and only if $\len _{\fp}
(\coker (\theta^{\natural} b))= \len _{\fp} (\coker (\theta^{\natural}
b'))$ for all $\fp$ in the inverse image under $\theta$ of each
singular point.

\begin{proposition}
Automorphisms of $(X,\ce,b)$ that are trivial on $X$ are of the form $
\gamma =(\zeta_{1},\zeta_{2},\dots,\zeta_{l}),$ where for each $i,$
$\zeta_{i}^{r}=1,$ and each $\zeta_{i}$ corresponds to a connected
component $ X^{\nu}_{i}$ of the curve $X^{\nu}.$

In other words, if $U_r$ is the group of \rth-roots of unity in $k,$
and $\Gamma(X^{\nu})$ is the dual graph of $X^{\nu},$
then $$\aut_X(\ce,b)= H^0(\Gamma(X^{\nu}), U_r).$$
\end{proposition}

\begin{pf}
$\theta : X^{\nu} \rightarrow X$ makes $ \ce \cong \theta_{*}
\theta^{\natural}
\ce,$ and $ \theta^{\natural}$ and $ \theta _{*}$
induce an equivalence of the categories of torsion-free rank-one
$\co_{X}$-{modules} which are singular at the double points normalized
by $\theta$ and invertible sheaves on $X^{\nu}.$ Therefore, it is
enough to study $
\theta^{\natural}(\ce).$ But automorphisms of line bundles on $X^{\nu}
$ are just given by $\nu$-tuples of $\zeta
\in k^{*}$; moreover, $\theta^{\natural}(\ce)^{r}=1$ implies that
$\zeta^{r}_{i}=1$ for all $i.$
\end{pf}

Note that, in general, isomorphisms of $(\ce,b)$ over $X$ must be
induced by isomorphisms of $\theta^{\natural} \ce,$ and therefore are
always constant on each connected component of $X^{\nu}$; namely, at a
singularity they are of the form $\Phi = \left(
\begin{array}{cc} \varphi_{+} & 0 \\ 0 & \varphi_{-} \end{array}
\right),$ with $\varphi_{+}$ and $\varphi_{-}$ in the base field $k.$

\subsection{Isomorphisms of Families of Spin Structures}

We can also say something about isomorphisms of spin structures in
general.  We have seen that over a field these are all constant, i.e.
at each singularity, any isomorphism of spin structures $$ \Phi : (E
(p,q), b)
\rightarrow (E (p',q'), b')$$ must be of the form $$\Phi = \left(
\begin{array}{cc} \varphi_+(0) & 0 \\ 0 & \varphi_-(0) \end{array}
\right).$$  Using the step-by-step method we can show that this is the
case over any complete local ring.  We just need to show that if
$\maxid I =0$ for some ideal $I$ in $R,$ and if $\Phi$ is constant
over $R/I,$ then $\Phi$ is constant over $R.$ But this follows because
$\varphi_+ = \varphi_+(0) + xi_+(x),$ for some $i_+(x) \in IR[[x]],$
and $\psi_+ = pi_+ =0;$ and similarly $\varphi_- = \varphi_-(0) +
yi_-(y),$ for some $i_-(y) \in IR[[y]]$ and $\psi_- =0.$ Therefore,
$$\Phi = \left(
\begin{array}{cc} \varphi_+(0) + xi_+(x) & 0
\\ 0 & \varphi_- (0) + yi_-(y) \end{array} \right),$$
and $b' = a'(x^u, p'x^{u-1}, \dots, y^v)$ is mapped to $b,$ namely $$b'
\circ \Phi^r = a'(x^u \varphi^r_+, x^{u-1} p' \varphi^{r-1}_+
\varphi_-, \dots, \varphi^r_- y^v) =b=a(x^u, x^{u-1}p, \dots, y^v).$$

So comparing the $b_0$ and $b_r$ terms,  $a' \varphi^r_+ = a =
a' \varphi^r_-,$ and $$\varphi^r_+ = (\varphi_+(0) + xi_+(x))^r= \varphi
_+(0)^r + rxi_+(x) \varphi_+(0)^{r-1}$$ $$ = \varphi^r_- = (\varphi_-(0)^r
+ yi_-(y))^r= \varphi_-(0)^r + ryi_-(y) \varphi_-(0)^{r-1}.$$

Equating terms with similar powers of $x$ and $y,$ and using the fact
that $\varphi_+(0)$ and $\varphi_-(0)$ are invertible, as is $r,$
we have  $i_+(x) = i_-(y) = 0,$ and $\Phi$ is constant.

\subsection{Automorphisms of Families of Quasi-Spin Structures}

Isomorphisms of quasi-spin curves are harder to classify than those of
spin curves, but if we limit ourselves to automorphisms, we can
completely classify these.

Given a quasi-spin structure $(\ce, b)$ on $\cx/S,$ with $S$ local and
complete, we want to study $\aut_{X}(\ce, b).$  First, we study the
local structure, namely, automorphisms of $(E(p,q), b)$ over $ A =
R[[x,y]]/pq-xy,$ with $b=(b_{0}, \dots, b_{r}).$  Here again, $R$ is a
complete local ring with maximal ideal $\maxid.$  Now, given $ \Phi
\in \aut_{A} (E,b) ,$ we know $\Phi = \left( \begin{array}{cc}
\varphi_{+} & \psi_{+}
\\ \psi_{-} & \varphi_{-} \end{array} \right )$ with $ b\circ \Phi^{r}
= b,$ and $$\varphi_{+}(0) p = \varphi_{-}(0) p \splice \varphi_{+}(0)
q = \varphi_{-} (0) q$$ $$ \varphi_{+} =
\varphi_{+}(0) + x \gamma_{+}, \psi_{+} = p \gamma_{+}, \splice
\gamma_+ \in R[[x]]$$ $$
\varphi_{-} = \varphi_{-}(0) + y \gamma_{-}, \psi_{-}=q \gamma_{-},
\splice \gamma_- \in R[[y]].$$

Now $b=(b_{0}, \dots, b_{r}) \equiv (\beta_{0} x^{i}, 0, 0, \dots, 0,
\beta_{r} y^{j})\bmod \maxid,$ $(i+j =r)$ and
$\beta_{r}, \beta_{0} \in
k^{*},$ and $\Phi \equiv \left( \begin{array}{cc}
\overline{\varphi}_{+} & 0
\\ 0 & \overline{\varphi}_{-} \end{array} \right ).$  By the results
of the previous section,  $\varphi^r_+ \equiv 1 \equiv
\varphi^r_-.$  In fact, this will hold for the whole family, i.e. we
can replace congruence with equality.

\begin{proposition}
$\Phi = \left( \begin{array}{cc}\varphi_{+} & 0
\\ 0 & \varphi_{-} \end{array} \right )$ with $\varphi^{r}_{+} =1,
\varphi_{-}^{r} =1.$
\end{proposition}
\begin{pf}
Using the step-by-step method we can assume the claim is true
$\bmod I $ for some ideal $I$ with $\maxid I=0.$ So
$\varphi^{r}_{+} = 1 +i,$ which implies that $(\varphi_+ - \frac{i}{r
\varphi^{r-1}_{+}})^r =1.$ So $\varphi_{+} = \zeta + i_{+}$ for some
$i_+$ in $I \cdot R[[x]],$ with $\zeta^r =1.$ Similarly, $\varphi_{-}
= \xi+i_{-}$ for some $i_-$ in $I \cdot R[[y]],$ with $\zeta^r =1.$
Thus $\Phi = \left( \begin{array}{cc}\varphi_{+} & 0
\\ 0 & \varphi_{-} \end{array} \right )$ and $b \circ \Phi^{r} = b$
implies that $$b_{0} = b_{0} (\zeta+i)^{r} = b_{0} (1+r
\zeta^{r-1}i_{+}), \splice b_{r} = b_{r} (1+r \xi^{r-1}i_{-}).$$ This
implies that $b_{0} r \zeta^{r-1} i_{+} = 0 = b_{r} r
\xi^{r-1}i_{-}.$ But $ b_{0} \equiv x^{i}
\overline{\beta}_{0} \bmod \maxid,$ with $ \overline{\beta}_{0} \in
k^{*},$ and since $i_+$ and $i_-$ annihilate $\maxid,$ this implies
that $b_{0} r \zeta^{r-1} i_{+} = x^{i} \beta_{0} (r \zeta^{r-1})
i_{+}$ for some $ \beta_{0} \in R^{\cross}$ lifting $\bar{\beta}_0,$
and similarly for $b_r.$ Since $1/r \in R,$ $i_{\pm} =0,$ and
$\varphi^r_+= \varphi^r_-=1.$ Moreover, $\varphi_+$ and $\varphi_-$
are in $R.$
\end{pf}

Now note also that $ \xi p = \zeta p,$ $\xi q = \zeta q.$ So $ p(\xi -
\zeta) = q(\xi - \zeta) = 0.$ But if $ \gamma := \xi - \zeta ,$ then
$(\xi + \gamma)^{r} = 1,$ which implies that $$r \gamma \xi^{r-1} +
\left( \begin{array}{c} r \\ 2
 \end{array} \right) \gamma^{2} \xi^{r-2} + \dots = 0.$$ And if
$(\gamma)$ is a proper ideal, then $\bmod (\gamma^{2})$ we get $r
\gamma \xi^{r-1} \equiv 0,$ i.e. if $1/r \in R,$ $ \gamma
\in (\gamma^{2}) ,$ which implies that $\gamma \in \bigcap_{n}
(\gamma^{n}) \subseteq \bigcap_{n} \maxid^{n} = 0.$  This implies that
$\gamma =0.$

So either
\begin{enumerate}
\item  $\gamma$ is invertible, hence $ p$ and $ q $ are zero, or
\item  $\gamma$  is zero and $ \varphi_+ = \varphi_-.$
\end{enumerate}

So at each singularity with at least one of $p $ and $q$ not zero,
$\aut (E(p,q),b) = U_{r} = \{\zeta \in R^{\cross} | \zeta^{r}
=1\}.$  And thus all automorphisms of $(\ce,b)$ are also in $U_{r}$ if
no singularities are of type $E(0,0).$ A singularity of type $(0,0)$
has automorphisms of type $(\xi, \zeta ) \in U_{r} \cross U_{r}.$
Normalizing $X$ at each singularity of type $(0,0)$ to get $X^{\nu}$
shows that $\aut (\ce, b)$ will be of type $(\xi_{1} ,
\dots, \xi_{m}),$ where $m$ is equal to the number of connected
components of $X^{\nu}.$

\subsection{Properties of the Isom Functor}\label{isom}

For any two quasi-coherent sheaves $\ce$ and ${\ce}'$ on a curve $X/B$
the functors $\Hom (\ce, {\ce}')$ and $\isom(\ce,{\ce}')$ are
representable (cf. \cite[7.7.8 and 7.7.9]{ega3} and \cite{lau}).  For
the $B$-scheme $V$ and map $ \Phi : \ce_{X_{V}}
\rightarrow {\ce '}_{X_{V}}$ on $ X_{V}$ which represent the functor
$\isom (\ce, \ce'),$ the condition that $ \Phi^{r}$ commutes with
$ b$ and $b'$ is clearly an open condition, and thus is representable
over $V.$ Moreover, the scheme representing the functor $T \mapsto
\Isom_{X_{T}}((\ce_{T},{b_{T}}), ({\ce '}_{T}, {b'}_{T})),$ is
an open subscheme of $
\Hom_{X_{T}}(\ce_{T}, {\ce'}_{T}),$ which is quasi-projective of
finite type.  Thus we have the following proposition.

\begin{proposition}
For any two quasi-spin structures $(\ce,b)$ and $({\ce}', b')$ over a
stable curve $X/B,$ the functor $T \mapsto
\isom_{X_{T}}((\ce_{T},{b_{T}}), ({\ce '}_{T}, {b'}_{T}))$ is
represented by a quasi-projective $B$-scheme of finite type.
\end{proposition}

The proof and proposition are also valid for isomorphisms of spin
structures and pure-spin structures.  Moreover, because the \isom
functor for stable curves over $S$ is representable by a quasi-projective
$S$-scheme of finite type, we actually have that for
any two (quasi/pure) spin curves $\frak S/T$ and $\frak S'/T'$ the
functor $T \mapsto \isom_{T\cross T'}(\frak S,\frak S')$ is also
representable by a quasi-projective $S$-scheme of finite type.

Not only is \isom\ representable, it is also unramified and finite,
as the next two propositions show.

\begin{proposition}
For any two quasi-spin curves (or spin curves or pure-spin curves)
$\frak{S}=(X,\ce,b)/T$ and $\frak{S}'=(X', {\ce}', b')/T',$ the scheme
$\Isom_{T \cross T'} (pr^*_1 \frak{S}, pr^*_2 \frak{S}')$ is
unramified over $T \cross T'.$
\end{proposition}

\begin{pf}
It suffices to show that for a ring $R$ with
square-zero ideal $I$ and for any two quasi-spin structures $(\ce,b)$
and $({\ce}',b')$ on a stable curve $X$ over $R$ with two isomorphisms
from $(\ce,b)$ to $({\ce}',b')$ which agree over $\bar{R} = R/I,$ the
two isomorphisms must then agree over $R.$  (We do not need to
consider isomorphisms of the underlying curve because the \isom\ functor
for stable curves is unramified.)  Since \Isom\ is a principal
homogeneous \Aut -space, we are reduced to showing that any automorphism
of $(\ce,b) $ which is the identity over $\bar{R}$ is the identity
over $R.$  But this follows easily from the fact that all
automorphisms of quasi-spin curves are constant and have \rth\ power
equal to the identity.
Therefore, \Isom\ is unramified.
\end{pf}

Since \Isom\ is of finite type and unramified, it is quasi-finite, so we
only need to check that it is proper to see that it is finite.

\begin{proposition}
For any two (quasi, pure) spin curves $\frak{S}=(X, \ce, b)/T$ and
$\frak{S'}=(X', {\ce}',b')/T',$ the scheme $\Isom_{T \cross T}
(pr_{1}^{*} \frak{S},$ $pr^{*}_{2} \frak{S}')$ is proper over $T
\cross T.$
\end{proposition}

\begin{pf}
We use the valuative criterion.  We must show that if we are given two
spin curves, quasi-spin curves, or pure-spin curves, $ \frak S=(X,
\ce, b)$ and $\frak{S}'=(X', \ce', b')$ both over $\spec R,$ where $R$
is a discrete valuation ring, and given an isomorphism $ \Phi_{\eta}
: {\frak S}_{\eta} \rightarrow {\frak S}'_{\eta}$ defined on the
generic fibres, then we can always extend $ \Phi_{\eta}$ to an
isomorphism $\Phi$ over all of $\spec R.$

We can also assume that $R$ is complete, and since for stable curves
the \isom\ functor is proper, we can assume that $X = X'.$  Now let
$Y$ be the $fpqc$ cover of $X$ given by $ Y = U \coprod \left(
\coprod_{\fp} \spec{\hat{\co}_{X,\frak{p}}}\right),$ with the union
being taken over all closed points $\fp$ of the singular locus of the
special fibre of $X,$ and $U$ the smooth locus of $X.$
If $\Phi_{\eta}$ extends to all of $Y,$ then it will in fact be
constant on all intersections $\spec{\hat{\co}_{X,\fp}}
\cross_{X} U,$ and these constant isomorphisms are uniquely determined by
$\Phi_{\eta},$ hence $\Phi_{Y}$ will descend to an extension of
$\Phi_{\eta}$ on $X.$

Thus we only need to consider the local situation; namely, about a
singular point of the special fibre.  This is the case where $$A = \rxy
\commasplice \ce = E(p,q) \splice \ce' = E (p', q')$$ with $pq = p'q'
= \pi.$ And we need to show that an isomorphism $\Phi_{\eta}$ on the
fibre over the field of fractions $K$ of $R$ extends to an isomorphism
on all of $A.$ $\Phi_{\eta}$ lifts to a map $ \tilde{\Phi}_{\eta}: (A
\tensor_{R} K) ^{\oplus 2} \rightarrow (A \tensor _{R} K) ^{\oplus 2},$
which induces the isomorphism $\Phi_{\eta}: E (p,q) \tensor_{R} K
\rightarrow E (p', q') \tensor_{R} K.$  Since $\Phi_{\eta} $ is
constant, $ \tilde{\Phi}_{\eta} $ is given as a matrix
$\tilde{\Phi}_{\eta} = \left ( \begin{array}{cc} \varphi_{+} & 0 \\ 0
& \varphi_{-} \end{array} \right) ,$ with $ \varphi_{\pm}
\in K.$  It suffices to show that $\varphi_{+}$ and $\varphi_{-}$ are
actually in $R.$  But to be an isomorphism,
$\tilde{\Phi}_{\eta}$ must be such that $ \tilde{b}' \circ
\tilde{\Phi}_{\eta}^{\tensor r} = \tilde{b}.$

And since $\tilde{b}= (b_{0}, b_{1}, \dots, b_{r})$ and $\tilde{b}'=
({b'_{0}},{b'_{1}}, \dots, {b'_{r}}),$ we have ${b'_{0}}
\varphi^{r}_{+} = b_{0}$ and $ {b'_{r}} \varphi_{-}^{r} = b_{r}.$ But
as we have seen, $b_{0}$ and ${b'}_{0}$ are both invertible in
$A[1/x],$ hence in $A[1/x]$ the constant $\varphi_{+}^{r} =
b_{0}/{b'}_{0} \in (\rxy) [1/x],$ and thus $\varphi_{+}^{r}
\in R,$ similarly for $\varphi_{-}^{r}.$ And $R$ is normal, hence the
$\varphi_{\pm}$ are in $R.$ And so $\Phi_{\eta}$ extends to all of
$\spec A,$ and thus to all of $X.$
\end{pf}

\section{Construction of the Stacks}

Fix $S$ to be a scheme of finite type over a field or over an excellent
Dedekind domain with $r$ invertible in $S.$  These conditions are
necessary for us to be able to use Faltings' theorem from
Section~\ref{falt} and to be able to use the standard theorems on
algebraic stacks (see \cite{dav}).

We have two main functors to consider, namely \qspin\ and \spinbar.
\qspin\ is the \'etale sheafification of the functor taking an
$S$-scheme $T$ to the set of isomorphism classes of quasi-spin curves
over $T.$  And \spinbar\ is the subfunctor of \qspin\ induced by
restricting to spin curves instead of quasi-spin curves.  Note that
for a quasi-spin structure the property of being a spin structure is
local on the curve in the \'etale topology; therefore, the property of
being a spin structure is independent of the choice of log structure.

We also will consider a third functor over \spinbar\ given by restricting
to pure-spin curves, namely those curves which, locally in the \'etale
topology, have the form $E(p,q)$ with $p = t^v,$ $q = t^u,$ for some
$t$ in the base, and $u+v=r.$  In particular, this means that $\pi$ is
$t^r.$  Up to \'etale covers this condition is also independent of the
particular choice of log structure and of the particular choice of $p$
and $q.$  We call this functor \purespin.  Of course, spin structures,
quasi-spin structures, and pure-spin structures are all the same thing
if the underlying curve is smooth.  And quasi-spin structures over a
reduced base (or if $\pi$ is not a zero divisor) are actually spin
structures.

The main result of this section is that \qspin, \spinbar, and
\purespin\ are all separated algebraic stacks, locally of finite type
over $\mgbar,$ the moduli space of stable curves, and \spin, the
moduli of smooth spin curves, is dense
in each of these.

The fact that \spin\ is dense in the stacks follows from
Proposition~{\ref{smoothdef}}.  To prove that these stacks are
algebraic, we need to do the following (see, for example, \cite[pp.
15--23]{dav}, or \cite{lau}):

\begin{enumerate}
\item Prove that the functors are limit preserving.
\item Provide a smooth cover $U$ of the stack.
\item  	Prove that for a fixed family of curves \cx\ over $T$
the functor $\isom_{U_{T} \cross U_{T}}(pr_{1}^*,pr_{2}^*)$ is
representable by a scheme (it is clearly a groupoid).
\item Prove that the stacks are separated by showing
that $\Isom_{U_{T} \cross U_{T}}(pr_{1}^*,pr_{2}^*)$ is actually proper
and finite over $U_T \cross U_T.$
\end{enumerate}

The last two conditions follow from the results of Section~\ref{isom}.
For the first two we begin by considering the stack \qspin.  Many
results on the other two stacks follow relatively easily from this
case.  The fact that the stack \qspin\ is limit preserving is a
straightforward consequence of the following theorem of Grothendieck
and the fact that the the condition on the length of the cokernel is
an open condition (c.f. Proposition~\ref{opencoker}), hence limit
preserving.

\begin{theorem}[{\cite[8.5.2]{ega4}}]
 Given a quasi-compact and quasi-separated scheme $S_0,$ and given a
 projective system $\{S_{\gamma} \}$ of $S_0$-schemes, relatively
 affine over $S_0,$ and quasi-coherent $\co_{S_{\gamma}}$-modules
 $\cf_{\gamma}$ and $\cg_{\gamma},$ with $ \cf_{\gamma}$ of finite
 presentation, the canonical homomorphism of groups $$
 \lim_{\rightarrow}(\Hom_{S_{\gamma}}(\cf_{\gamma},\cg_{\gamma}))
 \rightarrow \Hom_{S}(\cf,\cg)$$ is an isomorphism.  Here $S,$ $\cf,$
 and $\cg$ are the obvious limit objects.
\end{theorem}

All that remains is condition (2), i.e. to provide a smooth cover.

\subsection{A Smooth Cover of \qspin}

\begin{proposition}
Given a curve $X/B$ and an integer $N,$ sufficiently large, the
functor taking a $B$-scheme $T$ to the set of all triples
$(\ce,b,(e_1, \dots, e_n))$ where $(\ce,b)$ is a quasi-spin structure
on $X_T,$ and $(e_1, \dots, e_n)$ is a basis for the module $\Gamma
(X_T, \ce_T {\tensor}
\omega_T^{\tensor N})$ on
$X_T$ is representable.
\end{proposition}

\begin{pf}
Any quasi-spin structure $(\ce,b)$ on $X$ must have total degree $=
\frac{1}{r}(2g-2),$ and on its normalization $ \theta: \tilde{X}
 \rightarrow X$ $$ \theta^{\natural} \ce^{\tensor r} \cong \theta^*
\omega_{X/k} (-\sum (u_i {\fp}^+_i + v_i {\fp}^-_i)),$$ where the sum
is taken over all singularities $\{ \fp_i \}$ of $\ce$ and
$\{\frak{p}_i^+,
\frak{p}_i^-\}$ are the inverse images of $\frak{p}_i$ via $\theta.$
In particular, for any given irreducible component $X_j$ of $X,$ we
have $$\deg_{X_j} (\theta^{\natural}
\ce) \geq
\frac{1}{r} (\deg_{X_j} (\omega_{X/k}) -r \delta_j),$$ where
$\delta_j$ is the number of singularities of $X$ in $X_j.$ Since
$\deg_{X_j}(\omega_{X/k})$ is always positive, and since the total
number of singularities in a stable curve of genus $g$ is bounded by
$3g-3,$ we have $$\deg_{X_i} \theta^{\natural}
\ce \geq 1-\delta_i \geq 2 -3g.$$

Thus, given a very ample line bundle $\co(1)$ on $X,$ we can choose,
as in Theorem~\ref{thm:dsousa}, a fixed $m_0$ depending only on the
genus of $X$ such that $\ce \tensor \co(m) $ is generated by global
sections and has vanishing higher cohomology for all $m \geq m_0.$
Fix, once and for all, an integer $N$ large enough so that $\omega^N$
is very ample and $\ce \tensor \omega^N$ has all the desired
properties.  Now we can represent torsion-free rank-one sheaves with
bounded degree on each component by a subscheme of
$\quot_{\co_{X}^n/X/S}$ for some $n,$ i.e. there exists $ U_{1}
\hookrightarrow \quot_{\co_{X}^n/X/S}$ which represents the functor
$$ T \mapsto \{ \cf , ( e_{1} \dots e_{n})\},$$ where $ \cf$ is a
rank-one, torsion-free sheaf on $X_T$ with bounded degree on each
component, and $ (e_{1} \dots e_{n} )$ is a basis of $ \Gamma(X_{T},
\cf \tensor \omega^N)$ for $\omega^N$ sufficiently ample. So over
$X_1/U_1$ there is a universal pair $(\ce, ( e_{1} \dots e_{n})).$ And
to represent maps $ \ce^{\tensor r}
\rightarrow \omega _{X_{1}/U_{1}},$ take $$ V:= \Bbb V(\ce^{\tensor
r} \tensor \omega ^{-1}) := {\text{\bf Spec}}_{X_{1}} (
\text{Sym}_{\co_{X_{1}}}(\ce^{\tensor r} \tensor \omega ^{-1})),$$
So that $\Hom_{X_{1}}(Y,V) = \Hom _Y(\ce^{\tensor r}_Y, \omega_Y).$
So letting $V_T:= V \cross_{X_1} X_T = V
\cross_{U_1}T,$ we get that
 $$\Hom_{X_{1}}(X_{T},V)=
\Hom_{X_{T}}(X_{T},V_{T}) = (\prod_{X/S} V/X) (T)$$ is the functor we
want, and it is
representable because $X$ is flat and projective over $S$ (see \cite{fga}).

Now we have a universal triple $ \ce, ( e_{i} \dots e_{n}), b:
\ce^{\tensor r} \rightarrow
 \omega$ on $ X_{2}/U_{2}$ representing all maps $ \ce^{\tensor r}
\rightarrow \omega,$ and the additional condition that the cokernel of
$b$ is supported on the singular locus of $X_{2}$ is also
representable; namely, it is just the condition that $b$ is an
isomorphism on the complement of the discriminant locus, and this is
an open condition.  Finally, we need to represent the condition on the
cokernel, but this condition is open on the base, as proved in
Proposition~\ref{opencoker}.
\end{pf}

In general, for an arbitrary stable curve $\cx/T$ we have represented
by some scheme $U$ all quasi-spin structures $(\ce,b)$ on $\cx_U$ such
that $\ce \tensor \omega^N$ can be expressed as a quotient of
$\co^n_{\cx},$ together with a basis for the module $\Gamma(\cx, \ce
\tensor \omega^N).$  Moreover, at any closed point of $T$ all
quasi-spin structures on $\cx \cross_T \spec{\co_{T,t}}$ can be
expressed as such a quotient.  Therefore, at each point $u$ of $U$ the
complete local ring $\hat{\co}_{U,u}$ is a versal deformation of the
quasi-spin structure induced by $u.$  In particular, if the curve
$\cx$ we begin with is the universal curve over an \'etale cover $T$
of $\mgbar,$ the moduli stack of stable curves, then $U$ is a cover of
the stack \qspin.

\begin{proposition}
The scheme $U,$ which represents all quasi-spin structures $(\ce,b)$
on the universal curve $\cx/T$ together with a choice $(e_1, e_2,
\dots, e_n)$ of global sections which generate $\ce
\tensor \omega^N,$ is smooth over the stack of quasi-spin
structures on the universal curve.
\end{proposition}

\begin{pf}
We have to show that if an affine $T$-scheme $Y=\spec B$ has a
square-zero ideal $I \subseteq B$ and a quasi-spin structure $(\ce,b)$
on $X \cross_T Y$ such that $(\ce,b)$ restricted to $Y_0 = \spec{B/I}$
has a basis $(\bar{e}_1, \dots \bar{e}_n)$ for $\Gamma (X
\cross_T Y_0, \bar{\ce} \tensor
 \bar{\omega}^N),$ then $(\bar{e}_1, \dots, \bar{e}_n)$ lifts to a
basis of $\Gamma (X\cross_T Y, {\ce} \tensor {\omega}^N)$ on $Y.$
Namely, it suffices to show that if $pr: X \cross_T Y \rightarrow Y$
makes $pr_* (\bar{\ce} \tensor
\bar{\omega}^N)$ free on $Y_0,$ then the locally free sheaf $pr_* (\ce
\tensor \omega^N)$ is also free of the same rank as $pr_*(\bar{\ce}
\tensor \bar{\omega}^N).$  But this is clear because $\ce \tensor
\omega^N$ commutes with base change, and the exponential sequence $$ 0
\rightarrow \widetilde{M_n (I)}\rightarrow \cg l_n (\co_Y) \rightarrow
\cg l_n (\co_{Y_0}) \rightarrow 0 $$ shows that the kernel of the
homomorphism $H^1 (|Y_0|, \cg l_n (\co_Y)) \rightarrow H^1 (|Y_0|,
\cg l_n (\co_{Y_0}))$ is $H^1$ of the coherent sheaf
$\widetilde{M_n(I)}$ on an affine scheme, hence is zero.  Therefore,
any rank $n,$ locally free sheaf on $Y$ which restricts to a free
sheaf on $Y_0$ must be free already.
\end{pf}

This completes the proof that the stack \qspin\ is algebraic.

\subsection{\spinbar\ and \purespin}

To construct a smooth cover of \spinbar, we first take an arbitrary
smooth cover of \qspin, say $\frak{S}:U \rightarrow \qspin,$ together
with its canonical log structure.
We can also assume that $U$ is affine.  Thus all of our previous
descriptions of spin structures apply, and in particular, the
homomorphism $b: \ce^{\tensor r} \rightarrow \omega$ is given (\'etale
locally on $X/U$) as $b = (b_0, \dots, b_r).$  Also, we can still
write the $b_i$ as power series $$b_i = \sum_{n \geq 0} b_{i,n} x^n +
\sum_{m>0} b_{i,-m} y^m$$ with the same relations as before, and in
particular, up to suitable base extension and isomorphism of $E(p,q),$
$$p^u =q^v$$ and $$b_{0,i} = \pi^{u-i} b_{r,i-r} +
\sigma_i \text{ for } 0 < i \leq u \splice b_{r,i} = \pi^{i-u} b_{0,i} +
\sigma_i \text{ for } u \leq i < r.$$

Moreover, $b$ is a spin structure if and only if $\sigma_i = 0$ for all
$i.$  So the closed subscheme $V$ defined by the ideal generated by the
$\sigma_i$ actually represents the condition that $b$ is a
spin-structure.

It is clear that the spin curve $ \frak{S}_V=(\cx,\ce,b)_V$ over $V$
makes $V$ a cover of \spinbar.  To see that it is smooth over
\spinbar, note that for any spin curve $\frak{T}/T$ the scheme
$\Isom_{\qspin} (\frak{S}/U, \frak{T}/T)$ is isomorphic to
$\Isom_{\spinbar} (\frak{S}/U, \frak{T}/T),$ which is isomorphic to
$\Isom_{\spinbar} (\frak{S}_V/V, \frak{T}/T).$  Hence smoothness of
$U$ over \qspin\ implies smoothness of $V$ over \spinbar.

Alternately, we can consider, as in the construction of the versal
deformation of \qspin, a curve $\cx/T$ and a relatively torsion-free
sheaf $\ce,$ so that the pair $(\cx/T, \ce)$ is versal for stable
curves with rank-one, torsion-free sheaves with bounded degree on
each component.  Taking the canonical log structure and constructing
the canonical induced map $\ce \rightarrow \rho_* \cm,$ we can take the
scheme representing the property that $\cm$ is isomorphic to
$\omega_{X/T}$ to be our cover.
Since the property of being spin is independent of choice of log
structure, this is a cover of \spinbar.  Moreover, because it
represents all spin structure maps for $(\cx/T, \ce),$ it is smooth
over \spinbar.

The representability of \isom, as well as the other properties (finite
and unramified), all follow from the case of \qspin, hence \spinbar\
is an algebraic stack, locally of finite type over $S.$

Now to construct the stack \purespin, take, as above, the cover $V$ of
\spinbar\ together with its canonical log-structure on the universal
curve $\cx_V$ and isomorphisms $\ce \cong E(p,q)$ for $p,$ and $q$ in
$\co_V.$ The condition that $\frak{S}_V$ is pure is representable by
the relatively affine $V$-scheme $W:={\text{\bf Spec}}_V (\co_V
[\tau]/(p-\tau^v, q-\tau^u)).$ Again it is easy to verify that $W$ is
a smooth cover of \purespin, and that \purespin\ is algebraic.  We
have proved the following theorem.

\begin{theorem}
\qspin,  \spinbar,  and \purespin\ all form separated algebraic stacks
of finite type over \mgbar, and \spin\ is dense in each of these.
\end{theorem}

\subsection{Singularities and Smoothness of \spinbar\ and \purespin}

The deformation theory done in Section~\ref{deftheory} completely
describes the local structure of \spinbar, and in fact it shows that
\spinbar\ is relatively Gorenstein (i.e. \spinbar\ is Gorenstein if
the base $S$ is), namely it is enough to check the completion of the
stalks for an \'etale cover, (c.f.
\cite[Theorem 18.3]{matsumura}), and these are of the form
$$\hat{\co}_{S,s} [[P_1, Q_1, \dots, P_l, Q_l, t_{l+1}, \dots,
t_n]]/(P_i^{u_i}-Q_i^{v_i}).$$

Now it is enough to check the quotient $$\hat{\co}_{S,s}/(a_{1},
a_{2}, \dots, a_{m}),$$ where $\{ a_{i} \}$ are any
$\hat{\co}_{S,s}$-regular sequence.  And taking $a_{1}= P_{1},
a_{2}=P_{2}, \dots a_{l} = P_{l}, a_{l+1}= t_{l+1}, \dots, a_{n} =
t_{n},$ we are reduced to the case of $
\hat{\co}_{S,s}[[Q_{1}, Q_{2}, \dots Q_{l}]]/(Q_{i}^{v_{i}}).$
But this is Gorenstein because $\hat{\co}_{S,s}[[Q_{1}, \dots,
Q_{l}]]$ is, and the ring in question is just $\hat{\co}_{S,s}[[Q_{1},
\dots, Q_{l}]]$ modulo the regular sequence $ ( \{ Q_{i}^{v_{i}}\}).$

 \purespin\ provides a resolution of the singularities of \spinbar;
namely, \purespin\ is smooth over $S$ because the completion of any of
its local rings is of the form $$ \hat{\co}_{S,s} [[ \tau_1, \tau_2,
\dots, \tau_l, t_{l+1}, \dots, t_n]],$$ which is smooth over
$\hat{\co}_{S,s}.$

\section{Compactness}

The goal of this section is to prove the properness of the stacks
\qspin, \spinbar, and \purespin.  This is accomplished by studying the
boundary of these stacks, i.e. the degeneration of smooth spin curves
into spin structures on stable curves.  To prove that the stacks are
proper we will use the valuative criterion and the fact that smooth
spin curves are dense (c.f. Proposition~\ref{smoothdef}) to justify
checking the valuative criterion only in the case that the generic
fibre is smooth (c.f. \cite[pg. 109]{dm}).

\subsection{Extending Spin-Structures and Line Bundles}

Given a complete, discrete valuation ring $R$ with field of quotients
$K,$ and a $K$-valued point $\eta$ of $\qspin,$ corresponding to
$\frak S_{\eta} = (\cx_{\eta}, \ce_{\eta}, b_{\eta}),$ with
$\cx_{\eta}$ smooth over $K,$ we need to show that (up to finite
extension of $K$) there exists a quasi-spin curve $\frak S$ over $R,$
extending $\frak S_{\eta}.$ To this end, we construct a semi-stable
curve and line bundle which will give the desired extension when
contracted to its stable model, as in Section~\ref{geom}.

To begin, since $\mgbar$ is proper, there is a stable curve $\cx$
extending $\cx_{\eta}$ over $R.$ Take a uniformizing parameter $t \in
R$ and map $R$ to itself via $t \mapsto t^r.$ Pulling back $\cx$ along
this map yields another (singular) curve $\cx_r.$ Resolving the
singularities of $\cx_r$ by blowing up yields a semi-stable curve
$\tilde{\cx}$ with generic fibre $\cx_{\eta}$ (up to a finite
extension of $K$) and special fibre having chains of $n_i r-1$
exceptional curves over each singularity of $\cx_r.$ Here $n_i$ is the
order of the corresponding singularity of $\cx,$ namely $\cx$ has
local equation $ R[[x,y]]/xy-t^{n_i}.$

Now, since $\tilde{\cx}$ is regular, any line bundle on the generic
fibre will extend (but not uniquely) to the entire curve.  In
particular, there is some line bundle $\cl$ on $\tilde{\cx}$ which
extends $\cl_{\eta}.$ It is well-known that in such a case, any two
line bundles which agree on the generic fibre differ only by Cartier
divisors supported on the special fibre.  In other words, if
$\cm_{\eta} \cong \cn_{\eta}$ then $\cm \cong \cn \tensor \co (\sum
a_i X_i),$ where $X_i$ are the irreducible components of the special
fibre of $\tilde{\cx},$ and $a_i$ are integers.  In our case,
therefore, $\cl^{\tensor r} \cong \omega_{\tilde{\cx}}
\tensor \co (\sum a_i X_i)$ for some integers $a_i.$

Of course, if $\cl$ extends $\cl_{\eta}$ then any line bundle of the
form $\cl \tensor \co (\sum b_i X_i)$ also extends $\cl_{\eta}.$  The
following results show that there is a choice $\co (\sum b_i X_i)$ so
that $\cl':= \cl \tensor \co (\sum b_i X_i)$ is a line bundle with
degree zero on all but one exceptional curve per chain, has degree one
on the one remaining exceptional curve, and there exists an
$\co_{\tilde{\cx}}$-module homomorphism $\beta:{\cl'}^{\tensor
r}\rightarrow \omega_{\tilde{\cx}}$ which is an isomorphism everywhere
except on the exceptional curves where $\cl'$ has degree one.
Contracting all the exceptional curves of $\tilde{\cx}$ induces a spin
curve on $\cx_r,$ and hence an $R$-valued point of \qspin extending
$\frak{S}_{\eta}.$

\begin{proposition}
Given $\cl$ on $\tilde{\cx}$ such that $\cl^{\tensor r} \cong
\omega_{\tilde{\cx}} \tensor \co(\sum a_i X_i),$ the coefficients
$a_i$ which correspond to non-exceptional components of the special
fibre can all be assumed divisible by $r.$ In particular, the line
bundle $\displaystyle \cl':= \cl \tensor \co \left ( -\frac{1}{r}
\!\!\! \sum
\begin{Sb} X_i \text{not} \\
 \text{exceptional} \end{Sb} \!\!\! a_i X_i \right )$ has
${\cl'}^{\tensor r} \cong \omega \tensor \co (\sum e_j E_j)$ where all
the $E_j$ are exceptional curves.
\end{proposition}

\begin{pf}
Basic intersection theory shows that, for any curve $X_j,$ the degree
of $\co(\sum a_i X_i)$ on $X_j$ is $-a_j \delta_j + \sum a_i
\delta_{ij},$ where $\delta_j$ is the number of points in the
intersection of $X_j$ with the rest of the special fibre, and
$\delta_{ij}$ is the number of points in the intersection of $X_i$ and
$X_j.$ Now, on any given exceptional curve $E$ in a chain, with $E$
intersecting only two curves $C_1$ and $C_2,$ we have $\deg_E
(\omega(\sum a_i X_i)) = \deg_E(\co(\sum a_i X_i)) = -2e+c_1+c_2,$
where $e, c_1,$ and $c_2$ are the coefficients in the sum $\sum a_i
X_i$ of $E, C_1,$ and $C_2$ respectively.  Moreover,
$\deg_{X_i}(\omega(\sum a_i X_i)) = r \deg_{X_i} \cl \equiv 0
\pmod r$ for every $X_i.$ So, in particular, $c_1 + c_2 \equiv 2e
\pmod r.$  Now, given a chain of exceptional curves $E_1, \dots,
E_{nr-1},$ and the two non-exceptional curves $C$ and $D$ that the
chain joins, if their associated coefficients are $e_1, e_2, \dots,
e_{nr-1}, c, $ and $d,$ respectively, then we must have $e_2 \equiv
2e_1-c,$ $e_3 \equiv 2e_2-e_1\equiv 3e_1-2c$ and $e_i \equiv
ie_1-(i-1)c,$ so that $e_{nr-1} \equiv (nr-1)e_1-(nr-2)c$ and $d
\equiv nre_1 -(nr-1)c \equiv c.$ Therefore, since the special fibre is
connected, and since all of the non-exceptional curves are joined by
exceptional chains, all of the coefficients of the non-exceptional
curves are congruent to $c$ for some choice of $c.$ But since the
divisor $(\sum X_i)$ is trivial, we can assume that at least one of
the coefficients of a non-exceptional curve is zero, hence all of them
are congruent to zero $\pmod r,$ and thus $\cl':= \cl \tensor (-\sum
\!\!\! \begin{Sb} X_i \text{not} \\
\text{exceptional} \end{Sb}\!\!\! (a_i/r) X_i)$ is a line bundle extending
$\cl_{\eta}$ such that ${\cl'}^{\tensor r} \cong \omega(\sum e_i E_i)$
and the $E_i$ are all exceptional curves.
\end{pf}

\begin{proposition}
If $\cl$ is a line bundle on $\tilde{\cx}$ such that $\cl^{\tensor
r}\cong \omega(\sum e_iE_i),$ with all of the $E$'s exceptional curves
in the special fibre, then there is a choice of coefficients
$\{{e'}_i\}$ such that ${e'}_i \equiv e_i \pmod r$ for every $i,$ and
the degree of $\omega(\sum {e'}_i E_i)$ is zero on every exceptional
curve except perhaps one per chain, where it has degree $r.$  In
particular the bundle $$\cl ': = \cl \tensor \co(\sum
(\frac{{e'}_i -e_i}{r}) E_i)$$ has degree zero on every exceptional
curve except perhaps one per chain, where it has degree one.  And
${\cl'}^{\tensor r} \cong \omega(\sum {e'}_i E_i).$
\end{proposition}

\begin{pf}
Because $\cl^{\tensor r} \cong \omega(\sum e_i E_i),$ the degree of
$\co(\sum e_i E_i)$ on each $E_i$ must be congruent to zero $\pmod r,$
and so for any particular chain $E_1, \dots E_{nr-1}$ we have $e_2
\equiv 2e_1,$ $e_3 \equiv 3e_1,$ and $e_i \equiv ie_1$ for each $i.$
Choose $0 \geq {e'}_1 > -r$ with ${e'}_1 \equiv e_1 \pmod r,$ and let
${e'}_i = ie'_1$ for $1 \leq i \leq n(r+e'_1).$ Choose ${e'}_i = ie'_1
+r$ for $n(r+e'_1) +1 \leq i \leq nr-1.$ This gives ${e'}_i \equiv e_i
\pmod r$ for all $i,$ $\deg_{E_j} \co(\sum {e'}_i E_i) =0$ for all $j
\neq n(r+e'_1)$ and on $E_{nr+ne'_1},$ the degree is
$-2e'_{n(r+e'_1)}+e'_{n(r+e'_1)-1} +
e'_{n(r+e'_1)+1}=(-2(n(r+e'_1))+n(r+e'_1)-1+n(r+e'_1)+1)+r,$ which is
$r.$
\end{pf}

Note also that all of the ${e'}_i$ in the previous proposition were
negative, thus there is a canonical inclusion map $\omega(\sum {e'}_i
E_i)
\hookrightarrow \omega.$

\subsection{The Stacks are Proper}

By the results of the previous section, we have for any complete
discrete valuation ring $R$ with a smooth spin curve $\frak S_{\eta} =
(\cx_{\eta}, \cl_{\eta})$ over its generic point $\eta,$ an extension
of the spin curve to a curve/line-bundle pair $(\tilde{\cx},\cl)$ over
$R$ (up to finite extension of the field of fractions) with the
following special properties.
\begin{itemize}
\item $\tilde{\cx}$ is semi-stable.
\item For any chain of
 exceptional curves in the special fibre of $\tilde{\cx},$ the degree
of $\cl$ is zero on every exceptional curve in the chain except
perhaps one, where it has degree one.
\item $\cl^{\tensor r} \cong \omega_{\tilde{\cx}}
(\sum e_i E_i)$ with all of the ${e}_i$
negative.
\end{itemize}
This means there is a natural map $$ \omega(\sum {e}_i
E_i) \hookrightarrow \omega,$$ inducing a spin structure on the stable
model $\cx_r;$ namely, if
$\theta$ is the contraction $\tilde{\cx} \rightarrow \cx_r$ then
$\theta_* \cl$ is a rank-one torsion-free sheaf, and the map
$$\cl^{\tensor r} \irightarrow \omega(\sum e_i E_i)\hookrightarrow
\omega_{\tilde{\cx}}$$ induces a spin structure map $$(\theta_*
\cl)^{\tensor r} \rightarrow \omega_{\cx_r}.$$
Thus $\frak{S}_{\eta}$ extends to a spin structure $\frak{S}$ over all
of $R,$ and the valuative criterion holds.  We have proven the
following theorem.

\begin{theorem}
The stack \qspin is proper over $S.$
\end{theorem}

Since \spinbar is a closed subscheme of \qspin, and since it surjects
to \mgbar, it is also a compactification of \spin\ over \mgbar.  And
since $\purespin \rightarrow \mgbar$ is surjective, and \purespin\ is
proper over \spinbar, the stack \purespin\ is another compactification
of \spin\ over \mgbar.

\section*{Conclusion}

We have constructed three algebraic stacks which compactify the moduli
space of spin curves.  The stack of quasi-spin curves, which, in some
sense, is easiest to construct, is not as easy to describe as the
substack of spin curves, which has nice (Gorenstein) singularities.
And these singularities are resolved by the stack of pure-spin curves.

In the special case of $2$-spin curves, many of the difficulties
disappear.  In particular, all three compactifications coincide.
Moreover, this compactification of $2$-spin curves can be shown to
agree with those of Cornalba \cite{corn:theta} and Deligne
\cite{deligne:letter}.

\section*{Acknowledgements}

I am grateful to all the people who have helped with this paper, in
particular to Gerd Faltings, for introducing me to moduli problems in
general and this problem in particular, and for his help and direction
throughout.  Thanks also to B. Speiser, S.  Zhang, N. Katz, R. Smith,
S. Voronov, A. Vaintrob, J.-F.~Burnol, and J. H. Conway for helpful
discussions, and to P. Deligne for providing me with a copy of his
letter to Yuri Manin.

\nocite{altman:kleiman,altman:kleiman2}

\end{document}